\documentclass[a4paper]{article}

\usepackage{latexsym}
\input{amssym}
\newcommand{\al}{\alpha}
\newcommand{\ga}{\gamma}
\newcommand{\ep}{\epsilon}
\newcommand{\io}{\iota}\newcommand{\WD}{\wedge}
\newcommand{\la}{\lambda}\newcommand{\om}{\omega}

\newcommand{\Ga}{\Gamma}\newcommand{\De}{\Delta}
\newcommand{\Th}{\Theta}\newcommand{\La}{\Lambda}

\newcommand{\SO}{\mathbf{SO}}\newcommand{\GL}{\mathbf{GL}}

\newcommand{\Di}{\mathbf}\newcommand{\Cal}{\mathcal}
\newtheorem{THEO}{Theorem}[section]
\newtheorem{PR}{Proposition}[section]
\newtheorem{DEF}{Definition}[section]
\newtheorem{RE}{Remark}[section]
\newtheorem{CO}{Corollary}[section]

\newtheorem{EX}{Example}[section]
\newcommand{\mzz}[4]{\left(
\scriptsize\begin{array}{cc}#1&#2\\#3&#4\end{array}\right)}

\newcommand\C{{\mbox{$\Bbb C$}}}
\newcommand\R{{\mbox{$\Bbb R$}}}

\begin{document}
\begin{sloppypar}
\title{Twistorial constructions of\\ spacelike surfaces in Lorentzian
4-manifolds}
\author{Felipe Leitner}
\date{February 17, 1999}
\maketitle
\begin{abstract}
We investigate the twistor space and the Grassmannian fibre bundle of a Lorentzian 4-space with natural almost optical structures and its induced CR-structures. The twistor spaces of the Lorentzian space forms $\R^4_1$, $\Di{S}^4_1$ and $\Di{H}^4_1$ are explicitly discussed. The given twistor construction is applied to surface theory in Lorentzian 4-spaces. Immersed spacelike surfaces in a Lorentzian 4-space with special geometric properties like semi-umbilic surfaces and surfaces that have mean curvature of vanishing lengths correspond to holomorphic curves in the Lorentzian twistor space. For the Lorentzian space forms $\R^4_1$, $\Di{S}^4_1$ and $\Di{H}^4_1$ those surfaces are explicitly constructed and classified.
\end{abstract}

\section{Introduction}
The aim of this paper is the investigation of spacelike surfaces in Lorentzian 4-manifolds with the aid of a so-called twistor construction in 4-dimensional Lorentzian geometry. The idea of the draft, that is presented here, comes from Riemannian twistor theory and its well known application to the theory of surfaces in Riemannian four-spaces. We give a short description of the Riemannian construction. 

Consider an oriented, 4-dimensional Riemannian manifold $(M^4,g)$. The tangent space at any point of the manifold $M$ is isometric to the Euclidean 4-space $\R^4$. The set of complex structures on $\R^4$ can be identified with the homogenous space $\GL(4,\R)/\GL(2,\C)$. There are two kinds of orthogonal complex structures on $\R^4$; we set $J_o:=\mzz{0}{-1}{1}{0}\in\frak{gl}(2,\R)$:
$$A_+:=\{J\in\SO(4)|\ \ \exists\ A\in\Di{SO}(4)\ \mathrm{with}\ 
J=A^{-1}\mzz{J_o}{0}{0}{J_o}A\},$$
$$A_-:=\{J\in\SO(4)|\ \ \exists\ A\in\Di{SO}(4)\ \mathrm{with}\ 
J=A^{-1}\mzz{J_o}{0}{0}{-J_o}A\}.$$
Both sets $A_+$ and $A_-$ of orthogonal complex structures are naturally identified with the homogenous space $\SO(4)/\Di{U}(2)$, which is the 2-sphere $\Di{S}^2$. We choose the set $A_-$ and define the twistor space $\Cal{A}_-(M)$ of $M^4$ due to Atiyah/Hitchin/Singer [AHS78] as the associated fibre bundle $$\Cal{A}_-(M):=SO(M)\times_{\SO(4)}A_-=SO(M)\times_{\SO(4)}\Di{S}^2$$
over $M^4$ consisting of orthogonal complex structures on $TM$. $SO(M)$ denotes thereby the bundle of positive oriented orthonormal frames on $M$. The twistor space $\Cal{A}_-(M)$ admits two natural almost complex structures $J_1$ and $J_2$. The first almost complex structure $J_1$ on $\Cal{A}_-(M)$ is integrable if and only if the Riemannian 4-manifold $M^4$ is self-dual. The second complex structure $J_2$ is never integrable. 

Consider now a conformal immersion $f:N^2\hookrightarrow M^4$
of a Riemannian surface $N^2$ into the Riemannian 4-space $M^4$. The image $df(T_nN)$ of the tangent space at any point $n\in N$ is an oriented 2-plane in the Euclidean vector space $T_{f(n)}M$. Any oriented 2-plane $V^2$ in the oriented Euclidean four-space $\R^4$ gives naturally rise to an orthogonal complex structure on $\R^4$, which is the rotation around the angle $\frac{\pi}{2}$ in positive direction on $V$ and negative direction on $V^{\bot}$ and which is an element of the set $A_-$. Hence, the 2-plane $df(T_nN)$ in $T_{f(n)}M$ corresponds uniquely to an element $J_n$ in the fibre of the twistor space $\Cal{A}_-(M)$ over the point $f(n)$. This gives rise to a natural lift of the immersion $f$ to the twistor space $\Cal{A}_-(M)$:
$$\begin{array}{cccc}\ga_f:&N&\rightarrow&\Cal{A}_-(M)\ ,\\&n&\mapsto&J_n\end{array}\qquad\begin{array}{ccc}
&&\Cal{A}_-(M)\\&\ \ \qquad\stackrel{\ga_f}{\nearrow}&\downarrow\\
&f:N\ \rightarrow&M\end{array}.$$
J. Eells and S. Salamon studied this `Gauss' lift in [ES85] and proved the following relation for minimal surfaces.
\begin{THEO}
Conformally immersed minimal surfaces in a Riemannian 4-space $M^4$ correspond bijectively to non-vertical, $J_2$-holomorphic curves in the twistor space $\Cal{A}_-(M)$ over $M^4$.
\end{THEO}
In particular, the lift of a conformally immersed surface to the twistor space is horizontal iff the immersed surface is superminimal. This means that in the case of a self-dual Riemannian 4-space $M^4$, superminimal surfaces can be constructed by horizontal holomorphic curves in the twistor space, which is a complex manifold. Those constructions of superminimal surfaces have been done by Th. Friedrich in [Fri84] and especially for the hyperbolic four-space in  [Fri97]. By this way a well known result is obtained, which says that a superminimal surface in the Euclidean 4-space has locally the form
$$\begin{array}{ccl}\C&\rightarrow&\C\times\C\cong\R^4\ ,\\
z&\mapsto&(z,f(z))\end{array}$$
where $f$ is a holomorphic function. Using the twistor space $\Di{P}^3(\C)$ of the sphere $\Di{S}^4$, R. Bryant proved in [Bry82] the global result that any Riemannian surface admits a conformal, superminimal immersion into $\Di{S}^4$. 

The idea of twistor theory and its application to surface theory in Lorentzian geometry is as follows. The twistor space $\Cal{Z}(M^4_1)$ of an oriented Lorentzian 4-manifold $M^4_1$ is defined to be the bundle of null directions in the tangent space $TM$ over $M$ (comp. [Nur96]):
$$\Cal{Z}(M^4_1):=SO(M^4_1)\times_{\SO(3,1)}P,$$ where $P$ is the space of null directions in the Minkowski space $\R^4_1$. Instead of almost complex structures, the twistor space $\Cal{Z}(M^4_1)$ admits natural almost optical structures, which are related to CR-structures. 

Consider the conformal immersion $f$ of a Riemannian surface $N^2$ into an oriented Lorentzian 4-space $M^4_1$,
$$f:N^2\hookrightarrow M^4_1.$$
In any point $f(n)\in M^4_1$ of the immersed surface, there is an ordered pair of normal null directions on the tangent space $df(T_nN)$ to the surface in $M^4_1$. By choosing one of these normal null directions, we obtain a natural lift of the immersion $f$ to the bundle of null directions $\Cal{Z}(M^4_1)$,
$$\begin{array}{ccc}
&&\Cal{Z}(M^4_1)\\&\ \ \qquad\stackrel{\ga_f}{\nearrow}&\downarrow\\
&f:N^2\ \rightarrow&M^4_1\end{array}.$$
Similar as in 4-dimensional Riemannian geometry, spacelike surfaces in $M^4_1$ with special geometric properties can be characterized and constructed by holomorphic curves in the almost optical manifold $\Cal{Z}(M^4_1)$. 

The paper is organized in two parts. The first part,
section 2, is concerned with the twistor space and the Grassmannian fibre bundle of a Lorentzian 4-space and its natural almost optical structures. 
Special attention is given to underlying CR-manifolds (Theorem \ref{T2.3}). The twistor spaces of the Minkowski space $\R^4_1$, the pseudosphere $\Di{S}^4_1$ and the pseudohyperbolic space $\Di{H}^4_1$ are discussed. 

In the second part, section 3, we apply the Lorentzian twistor construction to surface theory. The second fundamental form $II$ of an isometrically immersed spacelike surface $(N^2,h)$ decomposes to
$$II=h\otimes H_++h\otimes H_-+L_++L_-,$$
where $H_+$ and $H_-$ are the lightlike parts of the mean curvature. The vanishing of components in this decomposition is related to the holomorphicity of the Gauss lifts of the immersion (Proposition \ref{P3.1}). We prove that non-vertical, holomorphic curves in the twistor space over a Lorentzian 4-space $M^4_1$ project to semi-umbilic ($L_+=0$ or $L_-=0$) surfaces in $M^4_1$ and to surfaces that have mean curvature of vanishing lenghts (Theorem \ref{T3.2}).
For the Lorentzian space forms $\R^4_1$, $\Di{S}^4_1$ and $\Di{H}^4_1$, we construct explicitly any semi-umbilic surface, which has in addition mean curvature of vanishing lenghts.

\section{Lorentzian twistor space and optical geometry}
\subsection{Optical geometry and CR-geometry}
In this section, we recall some basic facts about CR-geometry and optical geometry. Special attention is given to the relation between them. Optical geometry is introduced by A. Trautman in [Tra85] and is also treated in [Nur96].

We fix two notations.
Let $M^n_1:=(M^n,g)$ denote an $n$-dimensional Lorentzian $C^{\infty}$-manifold, where `Lorentzian' means that the smooth metric tensor $g$ is of signature $(n-1,1)$. For any oriented semi-Riemannian manifold $(M^n,g)$ of signature $(n-r,r)$, we denote by $SO(M)$ the $\SO(n-r,r)$-principal fibre bundle of oriented (pseudo)-orthonormal frames over $M$.

Consider a $2r$-dimensional $C^{\infty}$-distribution $\Cal{H}\subset TM^n$, $2r\leq n$, on an $n$-dimensional $C^{\infty}$-manifold $M$. A bundle morphism 
$$J:\Cal{H}\rightarrow\Cal{H},\quad J^2=-id,$$
is called an almost complex structure on the subbundle $\Cal{H}\subset TM$. 
In case that an almost complex structure $J$ is given on a distribution $\Cal{H}\subset TM^{2n-1}$ with codimension $1$ in $TM$, the pair $(\Cal{H},J)$ is called an almost CR-structure on the $(2n-1)$-dimensional $C^{\infty}$-manifold $M^{2n-1}$. An almost CR-structure $(\Cal{H},J)$ on $M^{2n-1}$ is called integrable or just a CR-structure on $M$ iff
$$[JX,Y]+[X,JY]\in\Ga(\Cal{H})\quad\forall X,Y\in\Ga(\Cal{H})\quad \mathrm{and}$$
$$J([X,JY]+[JX,Y])-[JX,JY]+[X,Y]=0\quad\forall X,Y\in\Ga(\Cal{H}),$$
i.e. the Nijenhuis tensor vanishes on $\Cal{H}$.
\begin{RE}
\begin{enumerate}
\item
An almost CR-structure $(\Cal{H},J)$ on a Riemannian manifold $(M^{2n-1},g)$ is orthogonal iff $J$ on $\Cal{H}$ is orthogonal with respect to $g$ on $\Cal{H}$.
\item
A $C^{\infty}$-map $f$ between two almost CR-manifolds $(M,\Cal{H},J)$ and
$(\tilde{M},\tilde{\Cal{H}},\tilde{J})$ is called a CR-map iff
$$df(\Cal{H})\subset\tilde{\Cal{H}},\quad df\circ J=\tilde{J}\circ df.$$
\end{enumerate}
\end{RE}
\begin{EX} \label{E2.1}
Natural CR-structure on the unit sphere bundle of an oriented Riemannian 3-manifold\\[2mm] 
Let $(N^3,h)$ be a $3$-dimensional, oriented Riemannian manifold.
The unit sphere bundle of $N^3$ is given by
\begin{eqnarray*}S^2(TN)&:=&SO(N)\times_{\SO(3)}\Di{S}^2=SO(N)\times_{\SO(3)}\SO(3)/\SO(2)\\ &\ =&SO(N)/\SO(2)\subset TN.\end{eqnarray*}
Let $\pi:SO(N)/\SO(2)\rightarrow N$ denote the natural projection.
The Levi-Civita connection of $(N^3,h)$ decomposes the tangent bundle $TS^2(TN)$ into a horizontal and a vertical part:
$$TS^2(TN)=T^VS^2(TN)\oplus T^HS^2(TN).$$
On $S^2(TN)$ exists a natural smooth distribution $\Cal{H}^{S^2(TN)}\subset TS^2(TN)$ of codimension $1$ given in a point $l\in S^2(TN)$ by
$$\Cal{H}^{S^2(TN)}_l=\pi_*^{-1}((\R l)^{\bot})=T^VS^2(TN)\oplus (T^HS^2(TN)\cap\pi_*^{-1}(\R l^{\bot})).$$
Let $J^{\Di{S}^2}$ denote the standard $\SO(3)$-invariant, complex structure on $\Di{S}^2$. It exists a natural complex structure on the distribution $\Cal{H}^{S^2(TN)}\subset TS^2(TN)$, pointwise defined by
$$J^{S^2(TN)}_l=\pi_*^{-1}\circ J_l\circ\pi_*+[e]^{-1}\circ J^{\Di{S}^2}\circ [e],$$
where $[e]$ denotes the identification of fibres with $\Di{S}^2$ by an orthonormal frame $e$ and $\pi_*^{-1}\circ J_l\circ\pi_*$ is the horizontal lift of the orthogonal complex structure $J_l$ on $(\R l)^{\bot}\subset T_{\pi(l)}N$, which is the rotation around the angle $\frac{\pi}{2}$ in positive direction. This CR-structure $(\Cal{H}^{S^2(TN)},J^
{S^2(TN)})$ on $S^2(TN)$ is always integrable. Furthermore, the unit sphere bundle together with its natural CR-structure is conformally invariant. 
\end{EX}

An almost optical structure $\Cal{O}=(\Cal{K},\Cal{L},J)$ on a $2n$-dimensional 
$C^{\infty}$-manifold $M$ consists of distributions $\Cal{K}$ and $\Cal{L}$, 
where $\Cal{K}\subset\Cal{L}\subset TM$ are subbundles, $\dim\Cal{K}=1$, $\dim\Cal{L}=2n-1$,
and an almost complex structure $J$ on the quotient bundle $\Cal{L}/\Cal{K}$.
An almost optical structure $\Cal{O}=(\Cal{K},\Cal{L},J)$ is by definition integrable iff 
\begin{enumerate}
\item[(A)]
$$ [\Ga(\Cal{K}),\Ga(\Cal{L})]\subset\Ga(\Cal{L}),\qquad
\phi_t^{k*}J=J\ \ \forall k\in\Ga(\mathcal{K}),$$
where $\phi_t^{k}$ denotes the flow of the field $k$. This means, locally 
in every point $m\in M$ the almost optical structure $\Cal{O}$ may be pushed down to an almost CR-structure $(\Cal{H}^m,J^m)$ on a locally induced quotient manifold $U_m/\!\!\sim^{\Cal{K}},$ where $U_m\subset M$ is a suitable neighborhood of $m\in M$ and any two points $u_1,u_2\in U_m$ are $\sim^{\Cal{K}}$-related iff both belong to the same integral curve of the distribution $\Cal{K}$, 
\item[(B)] these locally induced CR-structures $(\Cal{H}^m,J^m)$, $m\in M$, are integrable (comp. [Nur96]).
\end{enumerate}
\begin{RE}\label{R2.2}
\begin{enumerate}
\item
On a Lorentzian manifold $M^{2n}_1:=(M^{2n},g)$, the optical structure $\Cal{O}=(\Cal{K},\Cal{L},J)$ is called orthogonal if $\Cal{K}$ is a null subbundle, $\Cal{L}=\Cal{K}^{\bot}$ is the orthogonal complement and $J$ is orthogonal with respect to the Riemannian 
metric on $\Cal{L}/\Cal{K}$ induced by $g$.
\item
A $C^{\infty}$-map $f:M\rightarrow \tilde{M}$ between two almost optical manifolds $(M,\Cal{K},\Cal{L},J)$ and $(\tilde{M},\tilde{\Cal{K}},\tilde{\Cal{L}},\tilde{J})$ is called an optical map iff
$$ df(\Cal{K})\subset\tilde{\Cal{K}},\quad df(\Cal{L})\subset\tilde{\Cal{L}},\quad df\circ J=\tilde{J}\circ df.$$
\item
An almost CR-manifold $(N^{2n-1},\Cal{H},J)$ gives rise to a canonical almost optical structure on the manifold $M:=\R\times N$. This almost optical structure on $M$ is defined by
$$\Cal{K}^M:=T\R,\quad \Cal{L}^M:=T\R\oplus\Cal{H},\quad J^M\cong J:\ \Cal{L}^M/\Cal{K}^M\cong\Cal{H}\rightarrow\Cal{L}^M/\Cal{K}^M\cong\Cal{H}.$$
\item
If an almost optical structure $(\Cal{K},\Cal{L},J)$ on $M$ is integrable or just satisfies condition A, the locally induced almost CR-structure on $U_m/\!\!\sim^{\Cal{K}}$, $m\in M$, is given by
$$\Cal{H}^m:=\tilde{\pi}_*(\Cal{L}),\quad J^m:=\tilde{\pi}_*\circ J\circ\tilde{\pi}_*^{-1},$$
where $\tilde{\pi}:U_m\rightarrow U_m/\!\!\sim^{\Cal{K}}$ is the natural projection. The open neighborhood $U_m\subset M$ of $m\in M$ may be chosen in such a way that $(U_m,\Cal{K}|_{U_m},\Cal{L}|_{U_m},J|_{U_m})$ is optical diffeomorphic to $\R\times
 U_m/\!\!\sim^{\Cal{K}}$ with the natural almost optical structure given as above.
\item 
Let $N^{2n-1}$ be a submanifold of codimension $1$ in the almost optical manifold $(M^{2n},\Cal{K},\Cal{L},J)$ such that
$$TN\oplus\Cal{K}|_N=TM|_N.$$
It follows that $\Cal{H}^N:=TN\cap\Cal{L}|_N$ is a distribution in $TN$ with codimension $1$ and $\Cal{H}^N$ is naturally identified with $\Cal{L}/\Cal{K}|_N$. Therefore, $J$ induces an almost complex structure $J^N$ on $\Cal{H}^N$. The
pair $(\Cal{H}^N,J^N)$ is a naturally induced almost CR-structure on $N^{2n-1}$. We say that $(N,\Cal{H}^N,J^N)$ is an almost CR-submanifold of the almost optical manifold $(M,\Cal{K},\Cal{L},J)$.
\item
Let $(N^{2n-1},\Cal{H}^N,J^N)$ be an almost CR-submanifold of the almost optical manifold $(M^{2n},\Cal{K},\Cal{L},J)$. If $(\Cal{K},\Cal{L},J)$ is integrable or just satisfies condition A, the almost CR-structure $(\Cal{H}^N,J^N)$ in $m\in N\subset M$ is
locally equivalent to the naturally induced almost CR-structure $(\Cal{H}^m,J^m)$ on the quotient manifold $U_m/\!\!\sim^{\Cal{K}}$ in the point $\{m\}\in U_m/\!\!\sim^{\Cal{K}}$. In case that condition A is not satisfied, the almost CR-structure that is 
induced on a submanifold $$N^{2n-1}\subset(M^{2n},\Cal{K},\Cal{L},J),\quad TN\oplus\Cal{K}|_N=TM|_N,$$ depends on the imbedding of $N$ in $M$.
\item
Let $(\Cal{K},\Cal{L},J)$ be an almost orthogonal optical structure on a Lorentzian manifold $M^{2n}_1$. For any submanifold  $N^{2n-1}\subset M^{2n}_1$ of codimension 1, whose induced metric $g|_{TN}$ is positive definite (spacelike), it holds $$TN\cap\Cal
{K}|_N=\{0\}.$$ Therefore, the almost orthogonal optical structure $(\Cal{K},\Cal{L},J)$ on $M$ induces an almost orthogonal CR-structure $(\Cal{H}^N,J^N)$ on $N$.
\item
An $1$-dimensional distribution $\Cal{K}$ on a manifold $M$ is called regular if there is a smooth differentiable structure on the quotient set $M/\!\!\sim^{\Cal{K}}$ such that 
$\tilde{\pi}:M\rightarrow M/\!\!\sim^{\Cal{K}}$ is a $C^{\infty}$-submersion. If $(M,\Cal{K},\Cal{L},J)$ is an optical manifold, $M/\!\!\sim^{\Cal{K}}$ admits globally a natural CR-structure.
\end{enumerate}
\end{RE}

Until now, we have defined complex, CR- and optical structures. We want to consider mappings between them. 
\begin{DEF}Let $(Q,J^Q)$ denote an almost complex manifold, $(N,\Cal{H}^N,J^N)$ an almost CR-manifold and $(M,\Cal{K}^M,\Cal{L}^M,J^M)$ an almost optical manifold.
\begin{enumerate}
\item
A $C^{\infty}$-map $f:Q\rightarrow N$ 
is called holomorphic if 
$$df(TQ)\subset\Cal{H}^N\quad \mathrm{and}\quad df\circ J^Q=J^N\circ df.$$
\item
A $C^{\infty}$-map 
$g:Q\rightarrow M$ is called holomorphic if
$$dg(TQ)\subset\Cal{L}^M,\quad dg(TQ)\cap\Cal{K}={0},\quad \hat{\pi}\circ dg\circ J^Q=J^M\circ\hat{\pi}\circ dg,$$
where $\hat{\pi}:\Cal{L}\rightarrow\Cal{L}/\Cal{K}$ is the natural projection. 
\item 
A $C^{\infty}$-map 
$h:N\rightarrow M$ is called holomorphic if 
$$dh(\Cal{H}^N)\subset\Cal{L}^M,\quad dh(\Cal{H}^N)\cap\Cal{K}^M={0}\quad \mathrm{and}\quad \hat{\pi}\circ dh\circ J^N=J^M\circ\hat{\pi}\circ df.$$
\end{enumerate}
\end{DEF}
In the notations of the previous definition, we have obviously
\begin{PR} \label{P2.1}
\begin{enumerate}
\item
If the mappings $f:Q\rightarrow N$ and 
$h:N\rightarrow M$ are holomorphic, the map $g:=h\circ f:Q\rightarrow M$ 
is also holomorphic.
\item
If $f:Q\rightarrow M$ is holomorphic and $f(Q)\subset N\subset M$, where $N$ is an almost CR-submanifold of $M$, then the map $f:Q\rightarrow N$ is holomorphic. In case that $(\Cal{K}^M,\Cal{L}^M,J^M)$ on $M$ satisfies condition A and $f:Q\rightarrow M$ 
is holomorphic, at least locally the map $\tilde{\pi}\circ f$ is holomorphic, i.e. for an open subset $U\subset Q$ with $f(U)\subset U_m,\ m\in U_m\subset M$, the map 
$$\tilde{\pi}\circ f:U\subset Q\rightarrow U_m/\!\!\sim^{\Cal{K}}$$ is holomorphic.
\item
Let $(\Cal{K}^M,\Cal{L}^M,J^M)$ on $M$ be integrable and let $k\in\Ga(\Cal{K}^M)$ be a vector field on $M$. If the flow $\phi_t^k$ of $k$ is defined on $]-\ep,\ep[\times U$, where $\ep>0$ and $U\subset M$ is an open subset, and if $g:Q\rightarrow M$ is holomorphic, then the map
$$g_t:=\phi_t^k\circ g|_{g^{-1}(U)}:g^{-1}(U)\subset Q\rightarrow M$$ is holomorphic for any $t\in\ ]-\ep,\ep[$.
\end{enumerate}\end{PR}

\subsection{Twistor space of a Lorentzian manifold $M^4_1$}
The twistor space of an oriented Lorentzian $4$-manifold $M^4_1$ is defined to be the fibre bundle of null directions in the tangent space $TM^4_1$ (comp. [Nur96]). The twistor space admits natural almost optical structures. The integrability, the conformal invariance and the underlying CR-hypersurfaces of these optical structures are investigated in this section. We begin with the discussion of the fibre type of the twistor bundle.

Consider the Minkowski space $\R^4_1$. Let $\{u_i:i=1,\ldots,4\}$ denote the standard basis of $\R^4_1$ with $<u_1,u_1>^4_1=-1$ and let $P$ denote the set of null directions in $\R^4_1$. The space $P$ is a submanifold of the projective space $\Di{P}^3(\R)$ and is defined to be the fibre type of the twistor space of an oriented Lorentzian $4$-manifold. 
There are several characterizations of the fibre $P$. 
For the first, the fibre $P$ may be written as homogenous space. The Lorentzian group $\SO(3,1)$ acts transitively on $P$. Let $H_+$ and $H_-$ denote the isotropy groups of this action resp. in $\R(u_1+u_2)$ and $\R(u_1-u_2)$. We obtain 
$P\cong\SO(3,1)/H_+\cong\SO(3,1)/H_-.$  The elements
$$\{E_{ij}:\ i,j=1,\ldots,4,\ i<j\},\qquad E_{ij}:=\left(\begin{array}{cc}\cdot&-\ep_j\\
\ep_i&\cdot\end{array}\right),$$ form a basis of the Lie algebra $\frak{o}(3,1)$ of $\SO(3,1)$. The Lie algebras $\frak{h}_+$ and $\frak{h}_-$ of $H_+$ resp. $H_-$ are then given by
$$\frak{h}_+=Span\{E_{12},E_{13}+E_{23},E_{14}+E_{24},E_{34}\},\ \
\frak{h}_-=Span\{E_{12},E_{13}-E_{23},E_{14}-E_{24},E_{34}\}.$$
It is $$\frak{o}(3,1)=\frak{m}_+\oplus \frak{h}_+=\frak{m}_-\oplus\frak{h}_-,$$ 
$$\frak{m}_+=Span\{E_{13}-E_{23},E_{14}-E_{24}\},\quad \frak{m}_-=Span\{E_{13}+E_{23},E_{14}+E_{24}\}.$$ These decompositions of the Lie algebra $\frak{o}(3,1)$ are not reductive. 

The space $P$ of null directions in $\R^4_1$ is also naturally identified with the positive resp. negative projective spinor modul. For this let $\De_+\cong\C^2$ and $\De_-\cong\C^2$ denote the positive resp. negative spinor modul (comp. [Baum81]). To any spinor $\psi_{\pm}\in\De_{\pm}$ corresponds the $\R$-linear mapping
$$\begin{array}{cccc}\hat{\psi}_{\pm}:&\R^4_1&\rightarrow&\De_{\pm}\ ,\\ &x&\mapsto& x\cdot\psi_{\pm}\end{array}$$ where $x\cdot\psi_{\pm}$ is the Clifford product. The kernel $\ker\hat{\psi}_{\pm}$ of this mapping is a null direction in $\R^4_1$. Moreover, the mapping
$$\begin{array}{cccl}\io:&\Di{P}(\De_{\pm})&\rightarrow& P\ \ ,\\
&[\psi_{\pm}]&\mapsto&\ker\hat{\psi}_{\pm},\quad \psi_{\pm}\in[\psi_{\pm}]\end{array}$$ is a diffeomorphism. 
The natural complex structures $J^{\Di{P}(\De_{\pm})}$ on $P\cong\Di{P}(\De_{\pm})\cong\Di{P}^1(\C)$ are invariant by the $\SO(3,1)$-action on $P$ and are given on $\frak{m}_{\pm}\cong T_oP$ by 
$$J^{\Di{P}(\De_+)}(E_{13}-E_{23})=-(E_{14}-E_{24}),\quad 
J^{\Di{P}(\De_+)}(E_{14}-E_{24})=E_{13}-E_{23},$$
$$J^{\Di{P}(\De_-)}(E_{13}+E_{23})=-(E_{14}+E_{24}),\quad 
J^{\Di{P}(\De_-)}(E_{14}+E_{24})=E_{13}+E_{23}.$$ The identification of $\Di{P}(\De_+)$ and $\Di{P}(\De_-)$ via the space of null directions $P$ is anti-holomorphic.

Any unit timelike vector $T\in\R^4_1$, $g(T,T)=-1$, gives an identification of the space of null directions $P$ and the 2-sphere $S^2(T^{\bot})\cong\Di{S}^2$ in the orthogonal complement $T^{\bot}\cong\R^3$:
$$\begin{array}{ccccc}\io:&P&\leftrightarrow&S^2(T^{\bot})&.\\
&\R(T+S)&\mapsto&S\end{array}$$

The positive resp. negative twistor space of a $4$-dimensional, oriented Lorentzian manifold $M^4_1:=(M^4,g)$ is defined to be the positive resp. negative projective spinor bundle  $$\Cal{Z}_+(M^4_1):=P(S_+)=SO(M^4_1)\times_{\SO(3,1)}\Di{P}(\De_+),$$
$$\Cal{Z}_-(M^4_1):=P(S_-)=SO(M^4_1)\times_{\SO(3,1)}\Di{P}(\De_-).$$
Let $\pi:\Cal{Z}_{\pm}(M)\rightarrow M$ be the natural projection. We may these bundles also interpret as the bundle of null directions in $TM^4_1$,
$$\Cal{Z}_{\pm}(M)=SO(M)\times_{\SO(3,1)}P=SO(M)
\times_{\SO(3,1)}\SO(3,1)/H_{\pm}.$$
On the twistor space $\Cal{Z}_+(M)$ are given the natural almost optical structures $$\Cal{O}_+^+=(\Cal{K}_+,\Cal{L}_+,J^+_+)\quad\mathrm{and}\quad
\Cal{O}_+^-=(\Cal{K}_+,\Cal{L}_+,J^-_+).$$ They are defined pointwise as follows. Let $[\psi_+]\in\Cal{Z}_+(M)$ be an arbitrary point. For a suitable orthonormal basis $s=(s_1,\ldots,s_4)\in SO(M)$ we may write 
$$[\psi_+]=[s,\R(u_1+u_2)]\in\Cal{Z}_+(M).$$ To $[\psi_+]\in\Cal{Z}_+(M)$ corresponds the orthogonal optical structure
$$\Cal{O}^{[\psi_+]}=(K^{[\psi_+]},L^{[\psi_+]},J^{[\psi_+]})$$
in $T_{\pi([\psi_+])}M$ that is defined by
$$ K^{[\psi_+]}=\R(s_1+s_2)\subset T_{\pi([\psi_+])}M,\qquad L^{[\psi_+]}=Span\{s_1+s_2,s_3,s_4\}\subset T_{\pi([\psi_+])}M\quad \mathrm{and}$$
$$J^{[\psi_+]}(s_3+K^{[\psi_+]})=s_4+K^{[\psi_+]},\qquad
J^{[\psi_+]}(s_4+K^{[\psi_+]})=-s_3+K^{[\psi_+]}.$$
The almost optical structure $\Cal{O}^+_+=(\Cal{K}_+,\Cal{L}_+,J^+_+)$ is given in $T_{[\psi_+]}\Cal{Z}_+(M)$ by
$$\Cal{O}_{+[\psi_+]}^+=\pi_*^{-1}\circ\Cal{O}^{[\psi_+]}\circ\pi_*-[s]^{-1}\circ J^{\Di{P}(\De_+)}\circ[s],$$
where $\pi_*^{-1}\circ\Cal{O}^{[\psi_+]}\circ\pi_*$ denotes the horizontal lift of the optical structure $\Cal{O}^{[\psi_+]}$ in $T_{\pi([\psi_+])}M$ to $T_{[\psi_+]}^H\Cal{Z}_+(M)$. The almost optical structure $\Cal{O}^-_+$ is given in $T_{[\psi_+]}\Cal
{Z}_+(M)$ by
$$\Cal{O}_{+[\psi_+]}^-=\pi_*^{-1}\circ\Cal{O}^{[\psi_+]}\circ\pi_*+[s]^{-1}\circ J^{\Di{P}(\De_+)}\circ[s].$$
Analogously, we may define the almost optical structures $\Cal{O}_-^+=(\Cal{K}_-,\Cal{L}_-,J_-^+)$ and $\Cal{O}_-^-=(\Cal{K}_-,\Cal{L}_-,J_-^-)$ on $\Cal{Z}_-(M)$.
Under the identification
$$\Phi:\Cal{Z}_+(M)\cong SO(M)\times_{\SO(3,1)}P\cong\Cal{Z}_-(M)$$ of the positive and negative twistor space, the almost optical structures $\Cal{O}_+^+$ and $\Cal{O}_-^+$ are related by
$$\Phi^*(\Cal{K}_-,\Cal{L}_-,J_-^+)=(\Cal{K}_+,\Cal{L}_+,-J^+_+).\\[2mm]$$
\begin{THEO}\label{T2.1}(comp. [Nur96]) Let $M^4_1$ be an oriented, $4$-dimensional Lorentzian manifold and let $\Cal{Z}_+(M^4_1)$ and $\Cal{Z}_-(M^4_1)$ be its twistor spaces.
\begin{enumerate}
\item
The almost optical structures $\Cal{O}_+^+$ on $\Cal{Z}_+(M^4_1)$ and $\Cal{O}_-^+$ on $\Cal{Z}_-(M^4_1)$ are integrable iff $M^4_1$ is conformally flat.
\item
The almost optical structures $\Cal{O}_+^-$ and $\Cal{O}_-^-$ are never integrable.
\end{enumerate}
\end{THEO}
\begin{RE}
\begin{enumerate}
\item
The almost optical structure $\Cal{O}_+^+$ is integrable in 
$[\psi_+]=[s,\R(u_1+u_2)]\in\Cal{Z}_+(M)$ iff 
$\R(s_1+s_2)$ is a principal null direction (comp. [O'N93], [Nur96]). In terms of the Riemannian curvature tensor $R$, the integrability condition in this point is equivalent to
$$R_{1413}+R_{2413}+R_{1423}+R_{2423}=0,$$ $$R_{1414}+R_{2414}+R_{1424}+R_{2424}-R_{1313}-R_{2313}-R_{1323}-R_{2323}=0.$$
\item
The almost optical structures $\Cal{O}_+^-$ and $\Cal{O}^-_-$ are not only non-integrable, even more they don't induce locally any almost CR-structure. 
\end{enumerate}
\end{RE}

We consider the conformal invariance of the twistor spaces with its almost optical structures.
Let $\tilde{g}=\exp(2\rho)g$, where $\rho:M^4\rightarrow\R$ is a smooth function,  be a conformally equivalent metric to $g$ on $M^4$. The twistor spaces $\Cal{Z}_{\pm}(M,g)$ and $\Cal{Z}_{\pm}(M,\tilde{g})$ are naturally identified by
$$\begin{array}{cccc}\Th:&\Cal{Z}_{\pm}(M,g)&\cong&\Cal{Z}_{\pm}(M,\tilde{g})
\quad ,\\
&[s\cdot A,\R(u_1\pm u_2)]&\leftrightarrow&[\exp(-\rho)s\cdot A,\R(u_1\pm
u_2)].\end{array}$$
\begin{THEO}
Let $(M^4,g)$ be an oriented Lorentzian $4$-manifold and let $\tilde{g}=\exp(2\rho)g$ be a conformally equivalent metric to $g$. The identification $\Th$ of the twistor spaces
$(\Cal{Z}_{\pm}(M^4,g),\Cal{O}_{\pm}^+)$ and
$(\Cal{Z}_{\pm}(M^4,\tilde{g}),\Cal{O}_{\pm}^+)$ is an optical diffeomorphism.
\end{THEO} 
\textbf{Proof:} 
The distribution $\Cal{K}_+$ is the lightlike geodesic spray of the Lorentzian manifold $M^4_1$. The lightlike geodesic spray is conformally invariant. Hence, the optical flag
$$\Cal{K}_+\subset\Cal{L}_+\subset T\Cal{Z}_+(M)$$ is conformally invariant. 
The comparison of the almost complex structures $J^{+g}_{\pm}$ and $J^{+\tilde{g}}_{\pm}$ on the screen space $\Cal{L}_+/\Cal{K}_+$ shows that both are identical.\hfill$\Box$\\ \\
The almost optical structures $\Cal{O}_+^-$ and $\Cal{O}_-^-$ are not conformally invariant. 

We are now interested in the underlying CR-hypersurfaces of $(\Cal{Z}_{\pm}(M),\Cal{O}_{\pm}^+).$
Let $T\in\Ga(TM)$, $g(T,T)=-1$, be a timelike unit vector field on $M^4_1$. The choice of such a field $T$ is equivalent to a $\SO(3)$-reduction $SO_T(M)$ of the frame bundle $SO(M)$ over $M$, i.e. equivalent to an imbedding
$$\io:SO_T(M)\hookrightarrow SO(M)$$ of a $\SO(3)$-principal fibre bundle over $M$ in $SO(M)$. In particular, the twistor bundle may then be written as 
$$\Cal{Z}_+(M)=SO_T(M)\times_{\SO(3)}P\cong SO_T(M)\times_{\SO(3)}\Di{S}^2.$$
The twistor fibre $\pi^{-1}(m)$ over a point $m\in M^4_1$ is identified via $T_m\in T_mM$ with the sphere 
$$S^2(T_m^{\bot})\subset T_m^{\bot}\subset T_mM,$$ i.e. we can think of the twistor bundle as the bundle of 2-dimensional unit spheres, which are orthogonal to the timelike vector field $T\in\Ga(TM).$
 
Consider now an oriented spacelike hypersurface $N^3$ in $M^4_1$, i.e. the restriction $g|_{N}$ is positive definite. 
There is an unique timelike unit normal field $T$ on $N$ such that $(T,s_2,s_3,s_4)$ is positive oriented on $M^4_1$ if $(s_2,s_3,s_4)\in SO(N)$. The field $T$ induces a $\SO(3)$-reduction of the restricted frame bundle $SO(M)|_N$. Then, we have
\begin{eqnarray*}\Cal{Z}_+(M)|_N&=&SO(M)|_N\times_{\SO(3,1)}P\cong SO_T(M)|_N\times_{\SO(3)} S^2(T^{\bot})\\&\cong& SO(N)\times_{\SO(3)}\Di{S}^2,\end{eqnarray*}
which is the natural identification $\Psi$ of the restricted twistor space with the unit sphere bundle $S^2(TN)$ over $N$.
The restriction $\Cal{Z}_+(M)|_N$ of the twistor bundle to $N$ is a submanifold of codimension $1$ in $\Cal{Z}_+(M)$ and obviously, it holds
$$T(\Cal{Z}_+(M)|_N)\cap\Cal{K}_+=\{0\}.$$ From Remark \ref{R2.2}.5, it follows that the almost optical structure $\Cal{O}^+_+$ induces the almost CR-structure 
$(\Cal{H}_+,J_+)$ on  $\Cal{Z}_+(M)|_N$.
\begin{THEO} \label{T2.3}
Let $M^4_1$ be an oriented, $4$-dimensional Lorentzian manifold and let $N^3$ be an oriented, 3-dimensional spacelike submanifold of $M^4_1$. 
\begin{enumerate}
\item
The almost CR-structure $(\Cal{H}_+,J_+)$ on the restricted twistor bundle $\Cal{Z}_+(M)|_N$ is always integrable.
\item
The CR-structure $(\Cal{H}_+,J_+)$ on $\Cal{Z}_+(M)|_N$ and the natural CR-structure $(\Cal{H}^{S^2(TN)},J^{S^2(TN)})$ on the unit sphere bundle $S^2(TN)$ over the Riemannian $3$-manifold $N$ are equivalent under the natural identification $\Psi$
if and only if $N^3\subset M^4_1$ is a totally umbilic hypersurface.
\end{enumerate}
\end{THEO}
\textbf{Remark:} An analogous theorem holds for the negative twistor space $(\Cal{Z}_-(M),\Cal{O}_-^+)$.\\[2mm]
\textbf{Proof:} We need to prove the integrabilitiy of $(\Cal{H}_+,J_+)$.
The distributions $\Cal{H}_+$ and $\Cal{H}^{S^2(TN)}$ are identical under the identification $\Psi$:
$$\Psi_*(\Cal{H}_+)=\Cal{H}^{S^2(TN)}.$$
Let $s=(s_2,s_3,s_4):U\subset N\rightarrow SO(N)$ be an local orthonormal frame field on $N$. The local frame $s$ gives a trivialization of the unit sphere bundle:
$$S^2(TN)|_U\cong U\times\Di{SO}(3)/\Di{SO}(2).$$
The  horizontal lifts of $s_i$, $i=2,3,4$, with respect to the Levi-Civita connection of $N$ are given in this trivialization by
\begin{eqnarray*}s_i^H&=&s_i-\sum_{2\leq k<l\leq 4}g(\nabla_{s_i}s_k,s_l)E_{kl}\SO(2)\\
&=&s_i-g(\nabla_{s_i}s_2,s_3)E_{23}\SO(2)
-g(\nabla_{s_i}s_2,s_4)E_{24}\SO(2).\end{eqnarray*}
The almost complex structure $J_+$ on $\Cal{H}_+$ is given on $\Cal{H}^{S^2(TN)}$ in $TS^2(TN)$ via the identification $\Psi$ by
$$J_+(s_3^H)=s_4^{H}-2g(\nabla_{s_3}s_1,s_4)E_{23}\SO(2)+(g(\nabla_{s_3}s_1,s_3)-g(\nabla_{s_4}s_1,s_4))E_{24}\SO(2),$$
$$J_+(s_4^H)=-s_3^{H}+2g(\nabla_{s_3}s_1,s_4)E_{24}\SO(2)+(g(\nabla_{s_3}s_1,s_3)-g(\nabla_{s_4}s_1,s_4))E_{23}\SO(2),$$
$$J_+(E_{23}\SO(2))=E_{24}\SO(2),\qquad J_+(E_{24}\SO(2))=-E_{23}\SO(2).$$
An easy calculation shows that the Nijenhuis tensor of the almost CR-structure $(\Cal{H}_+,J_+)$ on $S^2(TN)\cong\Cal{Z}_+(M)|_N$ is defined and vanishes always. The formulas above also show that $J_+$ and $J^{S^2(TN)}$ on $\Cal{H}^{S^2(TN)}$ are identical iff
$$g(\nabla_{s_3}s_1,s_3)=g(\nabla_{s_4}s_1,s_4)\quad\mathrm{and}\quad
g(\nabla_{s_3}s_1,s_4)=0$$
for any orthonormal basis $(s_2,s_3,s_4)$ on $N$.
But this condition just means that the second fundamental form $II$ of $N$ in $M^4_1$ looks like
$$II=g\otimes H,$$
where $H$ is the mean curvature, i.e. the hypersurface $N$ in $M^4_1$ is totally umbilic.\hfill$\Box$\\[3mm]
The almost CR-structure on the restriction $\Cal{Z}_+(M)|_N$ that is induced by the almost optical structure $\Cal{O}^-_+$ on $\Cal{Z}_+(M)$ is never integrable.\\[3mm]
\textbf{Examples}\\[2mm]
A. The twistor space $\Cal{Z}_+(\R_1^4)$ of the Minkowski space $\R^4_1$ and the optical structure $\Cal{O}^+_+$\\[1mm]
Let $\R_1^4=(\R^4,<,>^4_1)$ be the flat Minkowski space. The twistor space $\Cal{Z}_+(\R_1^4)$ of $\R_1^4$ is given as
$$\begin{array}{ccccccc}
\R_1^4\times\Di{P}(\De_+)&\cong&\R_1^4\times P&\cong&\R_1^4\times\SO(3,1)/H_+&\cong&\R_1^4\times\Di{S}^2\\
(x,[\psi_+])&\leftrightarrow&(x,ker\hat{\psi}_+)=(x,\R(u_1+S))
&\leftrightarrow&(x,AH_+)&
\leftrightarrow&
(x,S)\\&&=(x,\R[A\cdot(u_1+u_2)])\end{array}$$
Let $\{u_i:i=1,\ldots,4\}$ be the standard basis in $T_x\R^4_1\cong\R^4_1$ and let $\{u_i^A\}$ be the orthonormal basis 
transformed by $A\in\SO(3,1)$, $(u_i^A)=(u_1,\ldots,u_4)\cdot A$. The optical structure $\Cal{O}_+^+$ on $\Cal{Z}_+(\R_1^4)=\R_1^4\times\SO(3,1)/H_+$ is given in $[\psi_+]=(x,AH_+)\in\Cal{Z}_+(\R_1^4)$ by
$$\Cal{K}_{[\psi_+]}=\R(u_1^A+u_2^A),\qquad\Cal{L}_{[\psi_+]}=Span\{u_1^A+u_2^A,u_3^A,u_4^A\},$$
$$J^+_+|_{T^V\Cal{Z}_+(\R_1^4)}\cong -J^{\Di{P}(\De_+)},\qquad
\begin{array}{c}J_+^+(u_3^A+\Cal{K}_{[\psi_+]})=u_4^A+\Cal{K}_{[\psi_+]}\\[1mm]
J_+^+(u_4^A+\Cal{K}_{[\psi_+]})=-u_3^A+\Cal{K}_{[\psi_+]}\end{array}$$
The twistor space $(\Cal{Z}_+(\R_1^4),\Cal{O}_+^+)$ is optical diffeomorphic to $\R\times(\R^3\times\Di{S}^2)$ with the optical structure
$$(T\R,T\R\oplus H^{S^2(T\R^3)},J^{S^2(T\R^3)}),$$
where $(\R^3\times\Di{S}^2,H^{S^2(T\R^3)},J^{S^2(T\R^3)})$ is the unit sphere bundle over the Euclidean space $\R^3$ with the natural induced CR-structure on it (comp. Example \ref{E2.1}).\\[3mm] 
B. The twistor space $\Cal{Z}_+(\Di{S}^4_1)$ of the pseudosphere $\Di{S}^4_1$ and the optical structure $\Cal{O}_+^+$\\[2mm] 
The hypersurface 
$$\Di{S}^4_1:=\{x\in\R^5_1:\ <x,x>^5_1=1\}$$
in $\R^5_1$ is with the induced metric of $\R^5_1$ an oriented Lorentzian 4-space of constant sectional curvature $1$ and is called the 4-dimensional pseudosphere (comp. [O'N83]). 
The pseudosphere $\Di{S}^4_1$ is a symmetric space:
$$\Di{S}^4_1\cong\SO(4,1)/\SO(3,1),$$
where 
$$\begin{array}{cccc}\io:&\SO(3,1)&\hookrightarrow&\SO(4,1)\\
&(A_{ij})&\mapsto&\left(\scriptsize\begin{array}{ccc}A_{11}&0&A_{1i}\\
0&1&0\\A_{i1}&0&A_{ij}\end{array}\right)\end{array}$$ is the imbedding of the isotropy group $\SO(3,1)$. 

The twistor space $\Cal{Z}_+(\Di{S}^4_1)$ is a homogenous space given as 
$$\Cal{Z}_+(\Di{S}^4_1)=\SO(4,1)/H_+$$
and is diffeomorphic to
$$\R\times S^2(T\Di{S}^3)\cong\R\times\SO(4)/\SO(2)\cong\R\times\Di{S}^3\times
\Di{S}^2,$$
where $S^2(T\Di{S}^3)\cong\SO(4)/\SO(2)$ is the unit sphere bundle over the sphere $\Di{S}^3\cong\SO(4)/\SO(3)$.
The optical structure $\Cal{O}_+^+=(\Cal{K}_+,\Cal{L}_+,J^+_+)$ on $\Cal{Z}_+(\Di{S}^4_1)$ is $\SO(4,1)$-equivariant and is given as follows.
For the Lie algebra $\frak{o}(4,1)=Span\{E_{ij}:1\leq i<j\leq 5\}$ of $\SO(4,1)$, we have the decomposition
$$\frak{o}(4,1)=\frak{m}_+\oplus \frak{h}_+\oplus\frak{b},$$
where $\frak{o}(3,1)=\frak{m}_+\oplus \frak{h}_+$ and $\frak{b}=Span\{E_{12},E_{23},E_{24},E_{25}\}$.
The subspace $\frak{b}$ is $Ad(H_+)$-invariant
and on $\frak{b}$ exists an $Ad(H_+)$-equivariant optical structure $(\frak{k},\frak{l},\underline{J})$ defined by
$$\frak{k}:=\R(-E_{12}+E_{23}),\qquad \frak{l}:=Span\{E_{24},E_{25}\}\oplus\frak{k},$$ 
$$\begin{array}{cccc}
\underline{J}:&\frak{l}/\frak{k}&\rightarrow&\frak{l}/
\frak{k}.\\
&E_{24}+\frak{k}&\mapsto&\ E_{25}+\frak{k}\\
&E_{25}+\frak{k}&\mapsto&-E_{24}+\frak{k}\end{array}$$
The twistor bundle $\Cal{Z}_+(\Di{S}^4_1)$ decomposes in a horizontal and a vertical part:
$$T\SO(4,1)/H_+=T^H\SO(4,1)/H_+\oplus T^V\SO(4,1)/H_+.$$
The horizontal bundle is given by
$$T^H\SO(4,1)/H_+=\SO(4,1)\times_{Ad(H_+)}\frak{b}.$$
The distributions $\Cal{K}_+$ and $\Cal{L}_+$ of the optical structure $\Cal{O}^+_+$ are defined by
$$\Cal{K}_+=\SO(4,1)\times_{Ad(H_+)}\frak{k},\qquad \Cal{L}_+=(T^V\SO(4,1)/H_+)\oplus(\SO(4,1)\times_{Ad(H_+)}\frak{l}).$$
The complex structure $J^+_+$ on $\Cal{L}_+/\Cal{K}_+$ is given by
$$J^+_+|_{T^V\SO(4,1)/H_+}\cong-J^{\Di{P}(\De_+)}$$ and
$$J_+^+([A,l+\frak{k}])=[A,\underline{J}(l+\frak{k})],\quad A\in\SO(4,1),\ l\in\frak{l},$$ on $T^H\SO(4,1)/H_+$.

The distribution $\Cal{K}_+$ is regular on $\Cal{Z}_+(\Di{S}^4_1)$ and the underlying CR-manifold of $(\Cal{Z}_+(\Di{S}^4_1),\Cal{O}_+^+)$ is equivalent to the unit sphere bundle of $\Di{S}^3$ with natural CR-structure. We describe this CR-structure on $\SO(4)/\SO(2)\cong S^2(T\Di{S}^3)$. Let $\{E_{ij}\}$ be the standard basis in $\frak{o}(4)$ and 
$$\begin{array}{cccc}\io:&\SO(2)&\hookrightarrow&\SO(4)\\
&\left(\scriptsize\begin{array}{cr}\cos t&-\sin t\\ \sin t&\cos t\end{array}\right)&\mapsto&\exp tE_{34}\end{array}$$
the imbedding of $\SO(2)$. Let
$$\frak{h}:=Span\{E_{13},E_{14},E_{23},E_{24}\}\subset\SO(4)$$
be a subspace of $\frak{o}(4)$. The subspace $\frak{h}$ is $Ad(\SO(2))$-invariant and the linear map
$$\underline{J}:\frak{h}\rightarrow\frak{h},\qquad\begin{array}{ll} E_{13}\mapsto E_{14},&E_{23}\mapsto E_{24},\\E_{14}\mapsto -E_{13},&E_{24}\mapsto-E_{23},\end{array}$$ is an $Ad(\SO(2))$-equivariant complex structure on $\frak{h}$. The canonical CR-structure on $\SO(4)/\SO(2)$ is given by
$$\Cal{H}^{S^2(TS^3)}=\SO(4)\times_{Ad(\SO(2))}\frak{h},$$
$$\begin{array}{cccc}J^{S^2(TS^3)}:&\Cal{H}^{S^2(TS^3)}&\rightarrow&\Cal{H}^{S^2(TS^3)}.\\&[A,h]&\mapsto&[A,\underline{J}h]\end{array}$$
C. The twistor space $\Cal{Z}_+(\Di{H}^4_1)$ of the pseudohyperbolic space $\Di{H}^4_1$ and the optical structure $\Cal{O}^+_+$\\[2mm]
The hypersurface 
$$\Di{H}^4_1:=\{x\in\R^5_2:\ <x,x>^5_2=-1\}$$
in $\R^5_2$ is with the induced metric of $\R^5_2$ an oriented Lorentzian 4-manifold of constant sectional curvature $-1$ and is called the 4-dimensional pseudohyperbolic space. 
The pseudohyperbolic space $\Di{H}^4_1$ is a symmetric space:
$$\Di{H}^4_1\cong\SO(3,2)/\SO(3,1),$$
where 
$$\begin{array}{cccc}\io:&\SO(3,1)&\hookrightarrow&\SO(3,2)\\
&A&\mapsto&\left(\scriptsize\begin{array}{cc}1&0\\0&A\end{array}\right)\end{array}$$ is the imbedding of the isotropy group. 

The twistor space $\Cal{Z}_+(\Di{H}^4_1)$ is a homogenous space given as 
$$\Cal{Z}_+(\Di{S}^4_1)=\SO(3,2)/H_+$$
and is diffeomorphic to
$$\Di{S}^1\times S^2(T\Di{H}^3)\cong\Di{S}^1\times\SO_o(3,1)/\SO(2)\cong\Di{S}^1\times\R^3
\times\Di{S}^2,$$
where $S^2(T\Di{H}^3)\cong\SO_o(3,1)/\SO(2)$ is the unit sphere bundle over the hyperbolic space $\Di{H}^3\cong\SO_o(3,1)/\SO(3)$.
The optical structure $\Cal{O}_+^+=(\Cal{K}_+,\Cal{L}_+,J^+_+)$ on $\Cal{Z}_+(\Di{H}^4_1)$ is $\SO(3,2)$-equivariant. We describe $\Cal{O}^+_+$ as follows.
The Lie algebra of $\SO(3,2)$ decomposes to
$$\frak{o}(3,2)=\frak{m}_+\oplus \frak{h}_+\oplus\frak{b},$$
where the subspace $\frak{b}=Span\{E_{12},E_{13},E_{14},E_{15}\}$
is $Ad(H_+)$-invariant. The subspaces
$$\frak{k}:=\R(E_{12}+E_{13}),\qquad \frak{l}:=Span\{E_{12}+E_{13,}E_{14},E_{15}\}$$ 
of $\frak{b}$ and the complex structure
$$\begin{array}{cccc}
\underline{J}:&\frak{l}/\frak{k}&\rightarrow&\frak{l}/
\frak{k}.\\
&E_{14}+\frak{k}&\mapsto&E_{15}+\frak{k}\\
&E_{15}+\frak{k}&\mapsto&-E_{14}+\frak{k}\end{array}$$
form together an $Ad(H_+)$-equivariant optical structure $(\frak{k},\frak{l},\underline{J})$ on $\frak{b}$. The optical structure $\Cal{O}^+_+$ on $\Cal{Z}_+(H^4_1)=\SO(3,2)/H_+$ is then given by
$$\Cal{K}_+:=\SO(3,2)\times_{Ad(H_+)}\frak{k},\qquad
\Cal{L}_+= (T^V\SO(3,2)/H_+)\oplus(\SO(3,2)\times_{Ad(H_+)}\frak{l})$$
and 
$$J^+_+:\Cal{L}_+/\Cal{K}_+\rightarrow\Cal{L_+}/\Cal{K}_+,\qquad
\begin{array}{c}J^+_+|_{T^V\SO(3,2)/H_+}\cong-J^{\Di{P}(\De_+)}\\[1mm]
J_+^+([A,l+\frak{k}])=[A,\underline{J}(l+\frak{k})]\end{array},$$
where $A\in\SO(3,2)$ and $l\in\frak{l}.$
The locally induced CR-structures of $\Cal{O}^+_+$ are equivalent to the natural CR-structure on the unit sphere bundle of $\Di{H}^3$. We describe this CR-structure on $\SO_o(3,1)/\SO(2)\cong S^2(T\Di{H}^3)$. Let
$$\begin{array}{cccc}\io:&\SO(2)&\hookrightarrow&\SO_o(3,1)\\
&A&\mapsto&\left(\scriptsize\begin{array}{cc}I&0\\0&A\end{array}\right)
\end{array}$$
be the imbedding of $\SO(2)$ in $\SO_o(3,1)$. The pair $(\frak{h},\underline{J})$ defined by
$$\frak{h}:=Span\{E_{13},E_{14},E_{23},E_{24}\}\subset\frak{o}(3,1)
\quad\mathrm{and}$$
$$\underline{J}:\frak{h}\rightarrow\frak{h},\qquad\begin{array}{ll} E_{13}\mapsto E_{14},&E_{14}\mapsto -E_{13},\\E_{23}\mapsto E_{24},&E_{24}\mapsto-E_{23},\end{array}$$ is an $Ad(\SO(2))$-equivariant CR-structure on $\frak{c}=Span\{E_{13},E_{14},E_{23},E_{24},E_{34}\}\subset\frak{o}(3,1)$. The CR-structure $(\Cal{H}^{S^2(T\Di{H}^3)},J^{S^2(T\Di{H}^3)})$ is then given by
$$\Cal{H}^{S^2(T\Di{H}^3)}=\SO_o(3,1)\times_{Ad(\SO(2))}\frak{h},$$
$$\begin{array}{cccc}J^{S^2(TH^3)}:&\Cal{H}^{S^2(TH^3)}&\rightarrow&\Cal{H}^{S^2(TH^3)}\\&[A,h]&\mapsto&[A,\underline{J}h]\end{array},\quad A\in\SO_o(3,1),\ h\in\frak{h}.$$

\subsection{Grassmannian fibre bundle of a Lorentzian manifold $M^4_1$}
The Grassmannian fibre bundle of an oriented Riemannian 4-manifold $M^4$ is the bundle of oriented planes in the tangent bundle $TM^4$. Each one of the positive and negative Riemannian twistor space may be seen as one `half' of the Grassmannian fibre bundle. In 4-dimensional Lorentzian geometry this point of view of the Grassmannian is also suitable up to an slight modification

Let $G(3,1)$ denote the set of oriented spacelike planes in the Minkowski space $\R^4_1$. We call $G(3,1)$ the Grassmannian of $\R^4_1$. The Grassmannian $G(3,1)$ can also be viewed as the set of decomposable 2-forms with lenght $1$:
$$G(3,1)\cong\{E\in\La^2\R^{4*}_1|<E,E>=1,\ E\WD E=0\}.$$ Furthermore, it is $G(3,1)\cong\SO(3,1)/T$ a symmetric space, where 
$$\begin{array}{cccl}\io:&T=\SO(1,1)\times\SO(2)&\hookrightarrow&\SO(3,1)\\
&(A,B)&\mapsto& \mzz{A}{0}{0}{B}\end{array}$$
is the isotropy subgroup. The Lie algebra $\frak{o}(3,1)$ decomposes to $\frak{o}(3,1)=\frak{m}\oplus\frak{t}$, where $\frak{m}=Span\{E_{13},E_{14},E_{23},E_{24}\}$ and $\frak{t}=\{E_{12},E_{34}\}$. The tangent bundle of $G(3,1)$ is the associated bundle
$$TG(3,1)=\SO(3,1)\times_{Ad(T)}\frak{m}.$$
Let $B(X,Y):=-\frac{1}{4}tr(adX\circ adY),\ X,Y\in\frak{o}(3,1)$, be the Killing form on $\frak{o}(3,1)$. The restriction $B|_{\frak{m}}$ of the Killing form to $\frak{m}$ is an $Ad(T)$-invariant scalar product on $\frak{m}$ and induces a $\SO(3,1)$-invariant metric $b^G$ of signature $(2,2)$ on the Grassmannian $G(3,1)$.
The Grassmannian $G(3,1)$ admits also the $\SO(3,1)$-equivariant complex structure $J_1,$ which is induced by the $Ad(T)$-equivariant complex structure
$$\begin{array}{cc}J^{\frak{m}}:\frak{m}\rightarrow\frak{m},\ &\ J^{\frak{m}}(E_{13})=E_{14},\ \ J^{\frak{m}}(E_{14})=-E_{13},\\ &\ J^{\frak{m}}(E_{23})=E_{24},\ \ J^{\frak{m}}(E_{24})=-E_{23}\end{array}$$ on $\frak{m}$. 
We denote $J_2:=-J_1$, which is a further complex structure on $G(3,1)$.
\begin{RE}
The Lie algebra $\frak{o}(4)\cong\La^2\R^4$ is semisimple.
Therefore, the Grassmannian of the Euclidean space $\R^4$ is isometric to the Riemannian product of the positive and the negative twistor space fibre:
$$G_2(4)\cong\Di{S}^2\times\Di{S}^2.$$ 
The Lie algebra $\frak{o}(3,1)\cong\La^2\R^4_1$ of the Lorentz group is simple. The Grassmannian $G(3,1)$ of $\R^4_1$ isn't the product of the positive and negative twistor space fibre $P$. It holds only
$$G(3,1)\cong P\times P\backslash{diag}\quad\subset P\times P.$$ It doesn't exist a natural $\SO(3,1)$-invariant metric on $P$. 
\end{RE}

Let $M^4_1$ be an oriented Lorentzian 4-manifold.
We call 
$$\Cal{G}(M):=SO(M)\times_{\SO(3,1)}G(3,1)=SO(M)/T$$
the Grassmannian fibre bundle of $M^4_1$.
The Grassmannian fibre bundle projects naturally to the twistor spaces. Let $e\in SO(M)$ and $s=e\cdot A\in SO(M),\ A\in\SO(3,1)$, be arbitrary orthonormal frames. The projections are given by
$$\begin{array}{crcl}\al_+:&\Cal{G}(M)\qquad&\rightarrow&\mathcal{Z}_+(M)\cong SO(M)\times_{\SO(3,1)}P\quad ,\\
&[e,AT]=Span\{s_3,s_4\}&\mapsto&[e,AH_+]=\mathbf{R}(s_1+s_2)\end{array}$$  
$$\begin{array}{crcl}\al_-:&\Cal{G}(M)\qquad&\rightarrow&\mathcal{Z}_-(M)\cong SO(M)\times_{\SO(3,1)}P\quad .\\
&[e,AT]=Span\{s_3,s_4\}&\mapsto&[e,AH_-]=\mathbf{R}(s_1-s_2)\end{array}$$
Moreover, we have the imbedding
$$\begin{array}{ccccc}\al:&\Cal{G}(M)&\hookrightarrow&\mathcal{Z}_+(M)
[\!\times\!]\mathcal{Z}_-(M)&,\\
&p&\mapsto&(\al_+(p),\al_-(p))&\end{array}$$
where $[\!\times\!]$ denotes the fiberwise product over $M$. 
In words, an oriented spacelike plane in a tangent space $T_mM^4_1$ is uniquely determined by its two ordered normal null directions in $T_mM^4_1$. 

On the Grassmannian fibre bundle exist the natural almost optical structures
$$\Cal{O}_G^+=(\Cal{K}_G,\Cal{L}_G,J^+_G)\quad \mathrm{and}\quad
\Cal{O}_G^-=(\Cal{K}_G,\Cal{L}_G,J^-_G).$$ They are given as follows.
Let $p\in\Cal{G}(M)$ be a point. In a suitable orthonormal basis $s=(s_1,s_2,s_3,s_4)$, we have $p=[s,eT]\in\Cal{G}(M).$ It is 
$$\Cal{O}^p=(\R(s_1+s_2),(s_1+s_2)^{\bot},J^p)\ ,$$
$$J^p(s_3+\R(s_1+s_2))=s_4+\R(s_1+s_2),\quad
J^p(s_4+\R(s_1+s_2))=-s_3+\R(s_1+s_2),$$
an optical structure on $T_{\pi(p)}M$. The almost optical structure $\Cal{O}^+_G$ is given in $T_p\Cal{G}(M)$ by
$$(\Cal{O}_{G}^+)_p=\pi_*^{-1}\circ\Cal{O}^p\circ\pi_*-[s]^{-1}\circ J_2\circ[s]$$
and $\Cal{O}_G^-$ is given in $T_p\Cal{G}(M)$ by
$$(\Cal{O}_{G}^-)_p=\pi_*^{-1}\circ\Cal{O}^p\circ\pi_*+[s]^{-1}\circ J_2\circ[s].$$
The projections
$$\al_+:(\Cal{G}(M),\Cal{O}_G^+)\rightarrow(\Cal{Z}_+(M),\Cal{O}_+^+),\qquad
\al_+:(\Cal{G}(M),\Cal{O}_G^-)\rightarrow(\Cal{Z}_+(M),\Cal{O}_+^-)$$ are optical maps.
\begin{THEO} Let $M^4_1$ be an oriented, 4-dimensional Lorentzian manifold and let $\Cal{G}(M)$ be its Grassmannian fibre bundle.
\begin{enumerate}
\item
The almost optical structure $\Cal{O}_G^-$ on $\Cal{G}(M)$ is never integrable.
\item
The almost optical structure $\Cal{O}_G^+$ on $\Cal{G}(M)$ is integrable in a point $p=[s,eT]\in\Cal{G}(M)$ if and only if
$$R_{1414}+R_{2414}-R_{1313}-R_{2313}=0,\quad
R_{1424}+R_{2424}-R_{1323}-R_{2323}=0,$$
$$2R_{1314}+R_{2314}+R_{2413}=0,\quad
2R_{2324}+R_{1423}+R_{1324}=0.$$
\item
The almost optical structure $\Cal{O}_G^+$ is integrable on $\Cal{G}(M)$ iff
$M$ has constant sectional curvature.
\end{enumerate}
\end{THEO}
\textbf{Sketch of the proof:} The integrability conditions in a point may be obtained as in the proof of Theorem \ref{T2.1}. Obviously, a space with constant sectional curvature satisfies these conditions in any point. On the other side, it can 
easily be seen that the Riemannian curvature tensor $R$ of a space, which hasn't constant sectional curvature, doesn't satisfy these conditions in some point\hfill$\Box$\\[2mm]
The almost optical structure $\Cal{O}^+_G$ induces an almost CR-structure on the restriction $\Cal{G}(M)|_N$ of the Grassmannian fibre bundle to an oriented spacelike hypersurface $N$. As in Theorem \ref{T2.3} this almost CR-structure is always integrable. 

To the end of this section, we define a parametrization of natural metrics on the Grassmannian fibre bundle $\Cal{G}(M)$. Let $0\neq\la\in\R$ be a parameter and $b^G$ the $\SO(3,1)$-invariant metric tensor on the Grassmannian $G(3,1)$. The tensor
$$g_{\la}:=\pi^*g+\la\cdot[s]^*(b^G),$$ where $\pi:\Cal{G}(M)\rightarrow M$ is the natural projection and $s\in SO(M)$ is an arbitrary frame, is for any $\la$ a metric of signature $(5,3)$ on $\Cal{G}(M)$.
  
\section{Surface theory and Lorentzian twistor construction}
\subsection{Immersed spacelike surfaces}
We study the second fundamental form of isometrically immersed surfaces with Riemannian metric in a Lorentzian 4-space.

Let $(N,h)$ be an oriented, 2-dimensional Riemannian manifold and let $dN$ be the volume form in the orientation of $N$. The Hodge operator $*$ is defined by
$\nu\WD *\xi=h(\nu,\xi)dN$. There is an unique orthogonal complex structure $J^N$ on $N$ such that $*\om=-\om\circ J^N$ for any 1-form $\om$ on $N$. It exist locally complex coordinates $(x_1,x_2)$ of the Riemannian surface $(N,J^N)$. For a suitable $C^{\infty}$-function $F$, we have
$$h=F\cdot(dx_1\circ dx_1+dx_2\circ dx_2).$$
The Bochner-Laplace operator on $C^{\infty}(N)$ is defined by $\De=*d*d$. For any $C^{\infty}$-function $f$ on $N$, we hace locally
$$\De f=\frac{f_{x_1x_1}+f_{x_2x_2}}{F}.$$
Let $M^4_1:=(M^4,g)$ be an oriented Lorentzian 4-space and let $f:N^2\hookrightarrow M^4_1$ be an isometric immersion. It is $f(N)\subset M^4_1$ a spacelike surface. A local Darboux frame is a map
$$e=(e_1,\ldots,e_4):U\subset N\rightarrow SO(M)|_N$$
such that $(e_3,e_4)$ is locally a positive oriented orthonormal frame in $TN$. Let $TN^{\bot}$ denote the normal bundle of the immersion $f$ over $N$. The metric $g$ on $TM$ induces a metric $h^{\bot}$ of signature $(1,1)$ on $TN^{\bot}$. The normal bundle decomposes into two null line bundles
$$TN^{\bot}=TN^{\bot}_+\oplus TN^{\bot}_-,$$
where these line bundles over $N$ are locally given by
$$TN_+^{\bot}=\R(e_1+e_2),\qquad TN_-^{\bot}=\R(e_1-e_2)$$
with respect to a Darboux frame $e$. We call $TN_+^{\bot}$ and $TN^{\bot}_-$ the bundles of positive resp. negative oriented normal null directions over $N$. The second fundamental form of the immersion $f$ is defined to be 
the normal part of the Levi-Civita connection $\nabla$ on $M^4_1$,
$$II(X,Y)=\Cal{N}\nabla_XY,\qquad X,Y\in\Ga(TN).$$
Let $h_{ij}^{\al}:=\ep_{\al}g(\nabla_{e_i}e_j,e_{\al})$ denote the components of $II$ with respect to a Darboux frame $e$. The mean curvature vector $H$ of the isometric immersion $f$ is given by
\begin{eqnarray*}H&=&\frac{1}{2}trII=\frac{1}{2}g^{ij}h_{ij}^{\al}e_{\al}
:=H_++H_-\\
&=&\frac{h_{33}^1+h_{44}^1+h_{33}^2+h_{44}^2}{4}(e_1+e_2)+
\frac{-h_{33}^1-h_{44}^1+h_{33}^2+h_{44}^2}{4}(e_2-e_1),\end{eqnarray*} 
where $H_-\in TN_-^{\bot},\ H_+\in TN_+^{\bot}$.
Let $t=\frac{-ie_3^*+e_4^*}{\sqrt{2}}\in\C\otimes T^*M$. It is
$$\begin{array}{crccl}II&=&h\otimes H&+&Re\left[\left(\frac{1}{2}(-h_{33}^1+h_{44}^1-h_{33}^2+h_{44}^2)+i(h_{34}^1+
h_{34}^2)\right)\cdot
t\circ t\right]\otimes(e_1+e_2)\\[1.5mm]
&&&+&Re\left[\left(\frac{1}{2}(h_{33}^1-h_{44}^1-h_{33}^2
+h_{44}^2)+i(-h_{34}^1+
h_{34}^2)\right)\cdot
t\circ t\right]\otimes(e_2-e_1)\\[1.5mm] 
&:=&h\otimes H_+&+&h\otimes H_-+L_++L_-\end{array}$$ with 
$L_-\in Sym^2(TN)\otimes TN_-^{\bot},\ L_+\in Sym^2(TN)\otimes TN_+^{\bot}$.
The second fundamental form $II$ and the given decomposition of $II$ is independently defined of the conformal class of the surface $N$.
\begin{DEF}
The isometric (conformal) immersion $f:(N,h)\rightarrow (M,g)$ is called
$$\mathrm{positive\ (negative)\ semi-stationary}\Leftrightarrow H_-=0\ (H_+=0),$$
$$\mathrm{positive\ (negative)\ semi-umbilic}\Leftrightarrow L_-=0\ (L_+=0),$$
$$+\mathrm{isotropic}\Leftrightarrow H_-=L_-=0\Leftrightarrow II(V,W)\in TN_+^{\bot},$$ 
$$-\mathrm{isotropic}\Leftrightarrow H_+=L_+=0\Leftrightarrow II(V,W)\in TN_-^{\bot},$$ 
$$\mathrm{stationary}\Leftrightarrow H=0,\quad\mathrm{totally\ umbilic}\Leftrightarrow II=g\otimes H.$$
\end{DEF}
\begin{RE}\label{R3.1}
Consider an isometric immersion $f:(N,h)\hookrightarrow(M,g)$. Let $\tilde{g}:=\exp(2\rho)g$ be a conformally equivalent metric to $g$ on $M^4$. Denote by $\tilde{II}$ the second fundamental form of the isometric immersion $f:(N,f^*\tilde{g})\rightarrow(M,
\tilde{g})$. The comparison of the covariant derivatives gives (see [Bes87])
$$\tilde{\nabla}_XY=\nabla_XY+d\rho(X)Y+d\rho(Y)X-g(X,Y)grad\rho$$ and this yields
$$\tilde{II}=II-g\otimes\Cal{N}grad\rho.$$
If the metric $h$ is positive definite, we have
$$\tilde{II}=f^*
\tilde{g}\otimes[\exp(-2\rho)\cdot(H-\Cal{N}grad\rho)]+L_++L_-.$$
This shows that the vanishing of the components $L_+,\ L_-$ is invariant under conformal change of the metric $g$ on $M^4$. In particular, the property of an immersion to be totally umbilic is conformally invariant, whereas the stationary condition isn't conformally invariant.
\end{RE}

\subsection{Gauss lift of an immersed surface}
We define now the Gauss lifts of a spacelike immersed surface to the twistor space and the Grassmannian fibre bundle. Geometric properties of an immersed surface are related to the holomorphicity of its Gauss lifts. Moreover, we describe the surfaces that have a harmonic Gauss lift to the Grassmannian.

Let $(N^2,J^N)$ be a Riemannian surface and $f:N^2\hookrightarrow M^4_1$ a conformal immersion into an oriented Lorentzian 4-manifold. It is
$$df(T_nN)\in G(T_{f(n)}M,g_{(f(n)})\cong G(3,1)$$
for any $n\in N$ an oriented Euclidean 2-plane in $T_{f(n)}M$, i.e. an element in the Grassmannian $\Cal{G}(M)$ over $M$. The smooth mapping
$$\begin{array}{cccc}v_f:&N&\rightarrow&\Cal{G(M)}\ ,\\
&n&\mapsto& df(T_nN)\end{array}\qquad
\begin{array}{ccc}
&&\Cal{G}(M)\\&\ \ \qquad\stackrel{v_f}{\nearrow}&\downarrow\\
&f:N\ \rightarrow&M\end{array}$$
is called the Gauss map or Gauss lift of the immersion $f$. We define also the (projected) Gauss maps into the twistor spaces over $M^4_1$:
$$\ga_{f\pm}:=\al_{\pm}\circ v_f:N\rightarrow \Cal{Z}_{\pm}(M)\ ,\qquad
\begin{array}{ccc}
&&\Cal{Z}_{\pm}(M)\\&\ \ \qquad\stackrel{\ga_{f\pm}}{\nearrow}&\downarrow\\
&f:N\ \rightarrow&M\end{array}.$$
The image of a point $n\in N$ under the mappings $\ga_{f\pm}$ is the positive 
resp. negative normal null direction on $df(T_nN)$ in $T_{f(n)}M^4_1$.
\begin{DEF} Let $f:N^2\rightarrow M^4_1$ be a conformal immersion of a Riemannian surface into an oriented Lorentzian 4-space. The immersion $f$ is called. 
$$\begin{array}{l}
\Cal{O}^+_G\mathrm{-holomorphic\ if}\ v_f:(N,J^N)\rightarrow(\Cal{G}(M),
\Cal{O}_G^+)\ \mathrm{is\ holomorphic,}\\[1.5mm]
\Cal{O}^-_G\mathrm{-holomorphic\ if}\ v_f:(N,J^N)\rightarrow(\Cal{G}(M),
\Cal{O}_G^-)\ \mathrm{is\ holomorphic},\\[1.5mm]
\Cal{O}_+^{\pm}\mathrm{-holomorphic\ if}\ \ga_{f+}:(N,J^N)\rightarrow(\Cal{Z}_{+}(M),
\Cal{O}_+^{\pm})\ \mathrm{is\ holomorphic},\\[1.5mm]
\Cal{O}_-^{\pm}\mathrm{-holomorphic\ if}\ \ga_{f-}:(N,J^N)\rightarrow(\Cal{Z}_{-}(M),
\Cal{O}_-^{\pm})\ \mathrm{is\ holomorphic}.
\end{array}$$
The  conformal immersion $f:N^2\rightarrow M^4_1$ is called  $\pm$horizontal if the lift $\ga_{f\pm}$ to $\Cal{Z}_{\pm}(M)$ is horizontal, that 
means $$d\ga_{f\pm}(T_nN)\subset T^H_{f(n)}\Cal{G}(M)\quad\forall n\in N.$$
\end{DEF}

Consider the differentials of the lifts $v_f$ and $\ga_{f\pm}$. The horizontal parts $dv_f^H$ and $d\ga_{f\pm}^H$ are simply the horizontal lifts of $df$. For the vertical part $dv_f^V$ holds
$$dv_f^V=-\om_{13}E_{13}T+\om_{23}E_{23}T-\om_{14}E_{14}T+\om_{24}E_{24}T,$$
and the vertical parts $d\ga_{f+}^V$ and $d\ga_{f-}^V$ are given by
$$d\ga_{f+}^V=\frac{-\om_{13}-\om_{23}}{2}(E_{13}-E_{23})H_+
+\frac{-\om_{14}-\om_{24}}{2}(E_{14}-E_{24})H_+,$$
$$d\ga_{f-}^V=\frac{-\om_{13}+\om_{23}}{2}(E_{13}+E_{23})H_-
+\frac{-\om_{14}+\om_{24}}{2}(E_{14}+E_{24})H_-.$$
\begin{PR}\label{P3.1}
Let $f:N^2\rightarrow M^4_1$ be a conformal immersion. The following relations hold:
\begin{enumerate}
\item
$f$ is $\Cal{O}^-_+$-holomorphic $\Leftrightarrow\ H_-=0$ $\Leftrightarrow$ $f$ is positive semi-stationary\\[1.5mm]
$f$ is $\Cal{O}^-_-$-holomorphic $\Leftrightarrow\ H_+=0$ $\Leftrightarrow$ $f$ is negative semi-stationary
\item
$f$ is $\Cal{O}_+^+$-holomorphic $\Leftrightarrow\ L_-=0$\\[1.5mm]
$f$ is $\Cal{O}_-^+$-holomorphic $\Leftrightarrow\ L_+=0$
\item
$f$ is $\Cal{O}_G^+$-holomorphic $\Leftrightarrow$ $f$ is totally umbilic\\[1.5mm]
$f$ is $\Cal{O}_G^-$-holomorphic $\Leftrightarrow$ $f$ is stationary $\Leftrightarrow$ $f$ is harmonic
\item
$f$ is $+$horizontal $\Leftrightarrow\ H_-=L_-=0\Leftrightarrow\ $
$f$ is $+$isotropic\\[1.5mm]
$f$ is $-$horizontal $\Leftrightarrow\ H_+=L_+=0\Leftrightarrow\ $
$f$ is $-$isotropic\\
\end{enumerate}
\end{PR}
\begin{RE} We may also say that $\pm$isotropic conformal immersions are super-semi-stationary. This is analogous to the notation for surfaces in a Riemannian 4-manifold, that have a horizontal Gauss lift to the Riemannian twistor space. Those surfaces satisfy a stronger condition then minimality and are called superminimal.\end{RE}

For characterizing an isometric immersion, whose corresponding Gauss lift is harmonic, we have the following
\begin{THEO}
Let $f:(N^2,h)\rightarrow(M^4,g)$ be an isometric (conformal) immersion of an oriented Riemannian 2-manifold into an oriented 4-dimensional Lorentzian space form, whose constant sectional curvature is $S$. For $\la S\neq 12$, the Gauss lift 
$v_f:(N,h)\rightarrow(\Cal{G}(M),g_{\la})$ is harmonic iff $H=0$. In case that $\la S=12$, $v_f:(N,h)\rightarrow(\Cal{G}(M),g_{\la})$ is harmonic iff $\nabla H=0$.
\end{THEO}
\textbf{Proof:}
On the frame bundle $SO(M)$, we have a natural parametrization of metrics,
$$g_{\la}=\pi^*g+\la B,\quad \la\in\R\backslash\{0\},$$
such that $\tilde{\pi}:(SO(M),g_{\la})\rightarrow(\Cal{G}(M),g_{\la})$ is a Riemannian submersion. Let $e=(e_1,\ldots,e_4):U\subset N\rightarrow SO(M)$ be a Darboux frame. In addition, we can choose the Darboux frame $e$ such that the frames $(e_3,e_4)$ on $N$ and $(e_1,e_2)$ on $TN^{\bot}$ are parallel in a point $n\in U$:  
$$\nabla^Ne_i(n)=0,\quad \nabla^{TN^{\bot}}e_{\al}(n)=0,\qquad i=3,4,\ \al=1,2.$$
For the Levi-Civita connection $\nabla^{\la}$ of $(SO(M),g_{\la})$, we have
$$\nabla^{\la}_{e_i^H}e_j^H=(\nabla_{e_i}e_j)^H+\frac{1}{2}Z[e_i^H,e_j^H],\quad
\nabla^{\la}_{E_{mn}}E_{kl}=\frac{1}{2}[E_{mn},E_{kl}],$$
$$\nabla^{\la}_{e_i^H}E_{mn}=\frac{\la}{2}\sum_{k=1}^{4}\ep_k\cdot g(E_{mn},[e_k^H,e_i^H])e_k^H=\nabla^{\la}_{E_{mn}}e_i^H,$$
where $e_i^H$ denotes the horizontal lift of $e_i,\ i=1,\ldots,4$, to $SO(M)$.
For the tension field $\tau_e$ of the local frame $e$ in $n\in N$ holds
\begin{eqnarray*}\tau_e&=\quad &tr(\nabla de)=\nabla^{\la}_{de(e_3)}de(e_3)-de(\nabla_{e_3}e_3)+
\nabla^{\la}_{de(e_4)}de(e_4)-de(\nabla_{e_4}e_4)\\
&=\quad & H^H+\sum_{k<l}\ep_k\ep_l\cdot \left(e_3(g(\nabla_{e_4}e_k,e_l))+e_4(g(\nabla_{e_4}e_k,e_l))\right)E_{kl}\\
&\quad -&\la\sum_{k<l}g(\nabla_{e_3}e_k,e_l)\cdot\sum_{m}\ep_mR_{m3kl}\cdot e_m^H
-\la\sum_{k<l}g(\nabla_{e_4}e_k,e_l)\cdot\sum_{m}\ep_mR_{m4kl}\cdot e_m^H\\&\quad + &terms\ in\ E_{12},E_{34}.\end{eqnarray*}
With the assumption that $(M^4,g)$ has constant sectional curvature we obtain for the tension field of $v_f$:
\begin{eqnarray*}\tau&=&tr(\nabla dv_f)=d\tilde{\pi}(\tau_e)=H^{H}-\frac{\la S}{12}H^{H}\\&&
+e_3(h_{33}^1+h_{44}^1)E_{13}+e_4(h_{33}^1+h_{44}^1)E_{14}-e_3(h_{33}^2+h_{44}^2)E_{23}-e_4(h_{33}^2+h_{44}^2)E_{24}.\end{eqnarray*}
Thus $v_f$ is harmonic iff $H=0$ or $\la S=12,\ \nabla H=0$.
\hfill$\Box$

\subsection{Twistorial construction of spacelike surfaces}
The Proposition \ref{P3.1} of the previous section relates geometric properties of an immersed spacelike surface to the holomorphicity of its Gauss lift. In the following, the reconstruction of semi-stationary and semi-umbilic surfaces in Lorentzian 4-spaces by holomorphic curves in the Lorentzian twistor space is established. We can give a whole classification of semi-umbilic surfaces in conformally flat Lorentzian 4-spaces. In particular, we classify every isotropic surface in the Lorentzian space forms $\R^4_1$, $\Di{S}^4_1$ and $\Di{H}^4_1$.

\begin{THEO}\label{T3.2}
Let $(N^2,J^N)$ be a Riemannian surface and let $M^4_1$ be an oriented Lorentzian 4-manifold.
\begin{enumerate}
\item
If 
$$\ga:(N,J^N)\rightarrow(\Cal{Z}_+(M),\Cal{O}_+^-)$$ is a holomorphic map into the positive twistor space over $M^4_1$ such that $d\ga\neq 0$ on $N$ and $\ga$ is non-vertical, that means $d\ga(T_nN)\not\subset T^V_{\ga(n)}\Cal{Z}_+(M)\ \forall n\in N$, then the projection 
$$f:=\pi\circ\ga:N^2\rightarrow M^4_1$$ is a conformal immersion, which 
is positive semi-stationary ($H_-=0$).
\item
There is a bijective correspondence between $\Cal{O}_+^-$-holomorphic, non-vertical curves in $\Cal{Z}_+(M)$ and conformally immersed positive semi-stationary surfaces in $M^4_1$.
\end{enumerate}\end{THEO}
\begin{RE}\label{R3.3}
\begin{enumerate}
\item
There is also a bijective correspondence between $\Cal{O}_+^+$-holomorphic, non-vertical curves in $\Cal{Z}_+(M)$ and conformally immersed surfaces in $M^4_1$ with $L_-=0$ (semi-umbilic). 
\item
Conformally immersed stationary surfaces in $M^4_1$ correspond bijectively to $\Cal{O}_G^-$-holomorphic, non-vertical curves in the Grassmannian $\Cal{G}(M)$ over $M$.
\item
To a conformal immersion
$$f:(N^2,h)\rightarrow(P^3,g)$$ of an oriented surface into an oriented Riemannian 3-manifold exists the Gauss lift to the unit sphere bundle $S^2(TP)$ over $P$, given by $$\begin{array}{cccc}\ga_f:&N&\rightarrow&S^2(TP),\\
&n&\mapsto&l_{f(n)}\end{array}$$
where $l_{f(n)}$ is the positive oriented unit normal vector to $df(T_nN)$ in $T_{f(n)}P$.  
There is a bijective correspondence between non-vertical, holomorphic curves in the CR-manifold
$$(S^2(TP),H^{S^2(TP)},J^{S^2(TP)})$$ 
and totally umbilic surfaces in $P$.

If $P$ is an oriented spacelike hypersurface of the oriented Lorentzian 4-space $M^4_1$, then the Gauss lifts of $$f:N^2\hookrightarrow P^3\subset M^4_1$$ to $S^2(TP)$ and to $\Cal{Z}_+(M)$ are identical under the identification $$\Psi:S^2(TP)\cong\Cal{Z}_+(M)|_P.
$$ In case that $P^3\subset M^4_1$ is a totally umbilic hypersurface, the Gauss lift of $f:N^2\hookrightarrow P^3\subset M^4_1$ to $(\Cal{Z}_+(M),\Cal{O}_+^+)$ is holomorphic iff the Gauss lift to $(S^2(TP),H^{S^2(TP)},J^{S^2(TP)})$ is holomorphic, i.e. iff $N^2\subset P^3$ is totally umbilic.
\end{enumerate}
\end{RE}
\textbf{Proof of the Theorem \ref{T3.2}:} Let
$$\ga:(N,J^N)\rightarrow(\Cal{Z}_+(M),\Cal{O}_+^-)$$
be a non-vertical, holomorphic curve, $d\ga\neq 0$ on $N$. Obviously, the projected map 
$f:=\pi\circ\ga$ is an immersion. It is $d\ga(TN)\subset\Cal{L}_+$  and it follows
$$\pi_*(d\ga(T_nN))\subset L^{\ga(n)}\quad\forall n\in N.$$
Moreover, it is 
$\hat{\pi}\circ d\ga\circ J^N=J^-_+\circ\hat{\pi}\circ d\ga$,
which means that
$$pr_{\ga(n)}\circ df\circ J^N_n=J^{\ga(n)}\circ pr_{\ga(n)}\circ df:T_nN\rightarrow L^{\ga(n)}/K^{\ga(n)}\quad\forall n\in N,$$
where $pr_{\ga(n)}:L^{\ga(n)}\rightarrow L^{\ga(n)}/K^{\ga(n)}$ is the natural projection. $J^{\ga(n)}$ is orthogonal with respect to $g_{f(n)}$ and therefore the immersion $f:N^2\hookrightarrow M^4_1$ is conformal. The Gauss lift $\ga_f$
of the immersion $f=\pi\circ\ga$ is equal to the original map $\ga$ and the theorem follows from Proposition \ref{P3.1}\hfill$\Box$\\[3mm]
The statements in Remark \ref{R3.3} can be proved in the same way as Theorem \ref{T3.2}. For the third point of Remark \ref{R3.3} remember Theorem \ref{T2.3}.
Parts of Proposition \ref{P3.1} and Theorem \ref{T3.2} are also proved in [Bob98]. 

The existence of a holomorphic curve in an integrable optical manifold gives rise to whole flow of holomorphic curves (Propositon \ref{P2.1}). This is the idea to
\begin{THEO}\label{T3.3}
Let $M^4_1$ be an oriented, conformally flat Lorentzian 4-manifold and let $k\in\Ga(\Cal{K}_+)$ be a vector field in the distribution $\Cal{K}_+$ on the twistor space $\Cal{Z}_+(M^4_1)$.
If $f:N^2\hookrightarrow M^4_1$ is a conformally immersed Riemannian surface with $L_-=0$ (semi-umbilic), then the map
$$f_t:=\pi\circ\phi_t^k\circ\ga_f:N^2\rightarrow M^4_1,$$ where $\phi^k_t$ denotes the flow of the field $k$,  is at least locally for small $t$'s a conformal immersion of a positive semi-umbilic ($L_-=0$) surface in $M^4_1$.
\end{THEO}
\textbf{Proof:}
If $f:N^2\hookrightarrow M^4_1$ is a conformally immersed surface with $L_-=0$
the Gauss lift
$$\ga_f:(N,J^N)\rightarrow(\Cal{Z}_+(M),\Cal{O}_+^+)$$ is holomorphic. Because $M^4_1$ is conformally flat, $\Cal{O}_+^+$ on $\Cal{Z}_+(M)$ is integrable and therefore
$$\phi_t^k\circ\ga_f:(N,J^N)\rightarrow(\Cal{Z}_+(M),\Cal{O}_+^+)$$ is locally and for small $t$'s a holomorphic map (Proposition \ref{P2.1})\hfill$\Box$\\[3mm]
Theorem \ref{T3.3} may be interpreted as follows. Let $N^2$ be any oriented spacelike surface in $M^4_1$. The normal bundle $TN^{\bot}$ over $N^2$ in $M^4_1$ decomposes to the bundle of positive and negative normal null directions on the surface $N^2$ in $M^4_1$. Any positive oriented normal null vector on $N^2$ defines a null geodesic that starts from the surface $N^2$. The distribution $\Cal{K}_+$ on $\Cal{Z}_+(M)$ is the lightlike geodesic spray of $M$. Hence, Theorem \ref{T3.3} says that any smooth deformation of a semi-umbilic ($L_-=0$) surface $N^2$ along its positive oriented normal null geodesics is also positive semi-umbilic. 

One way of finding holomorphic curves in an almost optical manifold is that of finding holomorphic curves in a CR-hypersurface of the almost optical manifold. For example, if $P^3\subset M^4_1$ is a spacelike hypersurface in a Lorentzian 
4-manifold, then the unit sphere bundle $S^2(TP)$ with CR-structure $(\Cal{H}_+,J_+)$ is a CR-hypersurface in $(\Cal{Z}_+(M),\Cal{O}_+^+)$ (Theorem \ref{T2.3}). A non-vertical holomorphic curve in $(S^2(TP),\Cal{H}_+,J_+)$ is then a non-vertical, holomorphic curve in $(\Cal{Z}_+(M),\Cal{O}_+^+)$ and projects to a semi-umbilic surface in $M^4_1$. This fact is described in Remark \ref{R3.3} for the special case, when $P^3\subset M^4_1$ is totally umbilic.
\begin{CO}\label{C3.1}
Let $M^4_1$ be conformally flat. If $P^3\subset M^4_1$ is totally umbilic and any complete lightlike geodesic intersects the hypersurface $P^3$, then any positive semi-umbilic surface $N^2\subset M^4_1$ is locally a deformation of a totally umbilic surface $\tilde{
N}^2\subset P^3$ along the positive oriented normal null geodesics on $\tilde{N}^2$.
\end{CO}
\textbf{Proof:}
For any $n\in N^2$ exists an open neighborhood $U\subset N^2$ and a null field $k\in\Ga(TM)$, which is tangential to the positive oriented normal null geodesics that start from $U$ and intersect the hypersurface $P^3$, such that $\phi_{t_1}^{k}$ is a diffeomorphism between $U$ and 
$$\tilde{N}^2:=\phi_{t_1}^k(U)\subset P^3$$ for a suitable $t_1$. From Theorem \ref{T3.3} and Remark \ref{R3.3}.3, it follows that $\tilde{N}^2$ is totally umbilic in $P^3$. The open set $U\subset N^2$ is the deformation of the totally umbilic surface $\tilde{N}^2$ in $P^3$\hfill$\Box$\\[1mm] 
In particular, any semi-umbilic surface in a conformally flat space $M^4_1$ is locally the deformation of a totally umbilic surface $\tilde{N}^2$ along the positive oriented normal null geodesics of $\tilde{N}^2$.\\[3mm]
\textbf{Examples}\\[2mm]
A. Positive semi-umbilic and +isotropic surfaces in the Minkowski space $\R^4_1$\\[1mm]
Consider the flat Minkowski space $\R_1^4=(\R^4,<,>^4_1)$. 
The complete, totally umbilic spacelike hypersurfaces in $\R^4_1$ are the Euclidean and hyperbolic 3-spaces. The Euclidean 3-spaces in $\R^4_1$ are totally geodesic. We choose the Euclidean 3-space
$$\begin{array}{rcc}P^3=\R^3&\hookrightarrow&\R^4_1\\
(y_1,y_2,y_3)&\mapsto&(0,y_1,y_2,y_3)\end{array}$$
in the Minkowski space $\R^4_1$. Any maximal lightlike geodesic in $\R^4_1$ intersects $P^3=\R^3$. 

The complete totally umbilic surfaces in the Euclidean 3-space $\R^3$ are the 2-spheres and the planes. The conformal immersion
$$\begin{array}{cccl}j:&\Di{S}^2&\hookrightarrow&\R^4_1\\
&(y_1,y_2,y_3)&\mapsto&(0,\theta y_1,\theta y_2,\theta y_3),\quad y_1^2+y_2^2+y_3^2=1,\ 0<\theta\in\R\end{array}$$
of the 2-sphere is totally umbilic. The null vector field 
$$\frac{\partial}{\partial x_1}+y_1\frac{\partial}{\partial x_2}+y_2\frac{\partial}{\partial x_3}+y_3\frac{\partial}{\partial x_4}\in\Ga(T\R^4_1|_{j(\Di{S}^2)})$$
is normal to $j(\Di{S}^2)\subset\R^4_1$ and positive oriented. The positive oriented normal null geodesics that start from the sphere $j(\Di{S}^2)$ are given by
$$\ga_{(a,b,c)}(t)=(t,(t+\theta)a,(t+\theta)b,(t+\theta)c),\quad t\in\R,\
a^2+b^2+c^2=1.$$
Any deformation of $j(\Di{S}^2)\subset\R^4_1$ along the positive oriented normal null geodesics has the form
$$\begin{array}{cccl}j_{\la}:&\Di{S}^2&\hookrightarrow&\R^4_1,\\
&(y_1,y_2,y_3)&\mapsto&(\la,(\la+\theta)y_1,(\la+\theta)y_2,(\la+\theta)y_3)
\end{array}$$
where $\la$ is a smooth function on $\Di{S}^2$. By Theorem \ref{T3.3}, it follows that these conformal immersions of $\Di{S}^2$ are complete, semi-umbilic ($L_-=0$) surfaces in $\R^4_1$.

The isometric immersion
$$\begin{array}{cccc}i:&(\R^2,<,>)&\hookrightarrow&(\R_1^4,<,>^4_1)\\
&(z_1,z_2)&\mapsto&(0,0,z_1,z_2)\end{array}$$
of the Euclidean plan is totally geodesic. The positive oriented normal null geodesics that start from the plane $i(\R^2)$ are given by
$$\ga_{(a,b)}(t)=(t,t,a,b),\quad t\in\R.$$
Any immersed surface that is the deformation of $i(\R^2)\subset\R^4_1$ along the positive oriented normal null geodesics on $i(\R^2)$ has the form
$$\begin{array}{cccc}i_{\la}:&\R^2&\rightarrow&\R^4_1,\\
&(z_1,z_2)&\mapsto&(\la(z_1,z_2),\la(z_1,z_2),z_1,z_2)\end{array}$$
where $\la$ is a smooth function on $\R^2$. The immersion $i_{\la}$ is isometric and $i_{\la}(\R^2)\subset\R^4_1$ is a complete, semi-umbilic
$(L_-=0)$ surface in $\R^4_1$ for any function $\la$ on $\R^2$. From Corollary \ref{C3.1}, it follows that any positive semi-umbilic surface in $\R^4_1$ is locally an immersion of the form $j_{\la}$ or $i_{\la}$ up to an isometry of the Minkowski space $\R^4_1$. The property that a surface is positive semi-umbilic is independent of the conformal class of the ambient Lorentzian 4-space. Therefore, the immersions $i_{\la}$ and $j_{\la}$ describe locally any positive semi-umbilic surface in a conformally flat Lorentzian 4-space.

In fact, the immersions of the form $i_{\la}$ into $\R^4_1$ are not only positive semi-umbilic, but even $+$isotropic. Its second fundamental form is calculated as
$$II=<,>\otimes\left(\ \De\la\ \right)\cdot\frac{e_1+e_2}{2}\ \ +\ L_+.$$
\begin{THEO}
Any complete +isotropic surface in the flat Minkowski space $\R^4_1$ is isometric to $\R^2$ and up to an isometry of $\R^4_1$ a +isotropic surface has (locally) the form 
$$\begin{array}{cccc}i_{\la}:&\R^2&\rightarrow&\R^4_1,\\
&(z_1,z_2)&\mapsto&(\la(z_1,z_2),\la(z_1,z_2),z_1,z_2)\end{array}$$
where $\la$ is a smooth function on $\R^2$. 
\end{THEO}
\textbf{Proof:} Any $+$isotropic surface is semi-umbilic $(L_-=0)$. The immersions of the form $j_{\la}$ are obviously never $+$isotropic. Hence, locally any $+$isotropic surface is of the form $i_{\la}$. We need to prove the global statement. The surface $i_{\la}(
\R^2)\subset\R^4_1$ is isometric to $\R^2$ for any  function $\la$ on $\R^2$. These surfaces are complete. Any other $+$isotropic surface isn't a deformation of a complete Euclidean plan $\R^2$ in $\R^4_1$ and therefore isn't complete\hfill$\Box$
\begin{RE}
\begin{enumerate}
\item
If we consider $\R^2$ as the complex numbers, we may write any complete $+$isotropic surface in $\R^4_1$ as
$$\begin{array}{ccc}\C&\rightarrow&\R^4_1,\\
z&\mapsto&(\la(z),\la(z),z)\cdot A^{\bot}+v\end{array}$$
where $A\in\Di{SO}(3,1)$ and $v\in\R^4_1$. This is similar to the case of superminimal surfaces in the Euclidean 4-space, which can be described as a graph of a holomorphic function $f$:
$$\begin{array}{ccl}\C&\rightarrow&\C\times\C\cong\R^4\quad .\\
z&\mapsto&(f(z),z)=(Re(f),Im(f),z_1,z_2)\end{array}$$
\item
If the function $\hat{\la}$ on $\R^2$ is harmonic, i.e. $\De\hat{\la}=0$, the surface $i_{\hat{\la}}(\R^2)$ in $\R^4_1$ is  complete, stationary and positive semi-umbilic. Any such surface has this form up to an isometry of the Minkowski space $\R^4_1$.\\
\end{enumerate}
\end{RE}
B. Positive semi-umbilic and $+$isotropic surfaces in the pseudosphere $\Di{S}^4_1$\\[2mm] 
Consider the pseudosphere $\Di{S}^4_1$. Let $A$ denote the subset 
$$\{x\in\Di{S}^4_1\subset\R^5_1:\ x_3=1\}$$
of $\Di{S}^4_1$. The conformal diffeomorphism
$$\begin{array}{cccc}\Xi:&\Di{S}^4_1\backslash A&\leftrightarrow&\R^4_1\\
&(t,a,b,c,d)&\mapsto&\frac{1}{1-b}(t,a,0,c,d)\\
&\frac{2}{1+\|x\|^2}(x_1,x_2,\frac{\|x\|^2-1}{2},x_4,x_5)&\leftarrow&x=(x_1,x_2,0,x_4,x_5)\end{array}$$
is the stereographic projection of the pseudosphere $\Di{S}^4_1$ in the point
$(0,0,1,0,0)\in\Di{S}^4_1$ to the Minkowski space $\R^4_1$. The induced metric of constant sectional curvature 1 on $\R^4$ is
$$\frac{4}{(1+\|x\|^2)^2}\cdot <,>^4_1=\exp(2\rho)\cdot <,>^4_1,$$
where $\rho=\ln\frac{2}{1+\|x\|^2}$.
The image of the isometric imbedding
$$\begin{array}{clcl}\io:&\Di{S}^3\subset\R^4&\hookrightarrow&\Di{S}^4_1
\quad\subset\R^5_1\\
&y&\mapsto&(0,y)\end{array}$$
is a complete, totally geodesic spacelike hypersurface in $\Di{S}^4_1$. Any maximal lightlike geodesic in $\Di{S}^4_1$ intersects this hyperpshere. The complete totally umbilic surfaces in $\Di{S}^3$ are the 2-spheres
$$\begin{array}{cccc}i^c:&\Di{S}^2&\hookrightarrow&\Di{S}^3\subset\R^4,\\
&(y_1,y_2,y_3)&\mapsto&(c,\sqrt{1-c^2}\cdot y_1,\sqrt{1-c^2}\cdot y_2,\sqrt{1-c^2}\cdot y_3)\cdot A^{\bot}\end{array}$$
where $|c|\leq 1$ and $A\in\SO(4)$. The surface $i(\Di{S}^2):=i^0(\Di{S}^2)\subset\Di{S}^3$ is totally geodesic. The conformal immersion
$$\begin{array}{cccc}i_{\la}^c:&\Di{S}^2&\hookrightarrow&\Di{S}^4_1\\
&y=(y_1,y_2,y_3)&\mapsto&(\la,\la\sqrt{1-c^2}+c,(\sqrt{1-c^2}-\la c)y)\end{array}$$
is for any smooth function $\la$ on $\Di{S}^2$ a deformation of $i^c(\Di{S}^2)\subset\Di{S}^4_1$ along the positive oriented normal null geodesics. By Theorem \ref{T3.3} and Corollary \ref{C3.1}, it follows that any positive semi-umbilic surface in $\Di{S}^4_1$ is locally up to an isometry of $\Di{S}^4_1$ an immersion of the form $i_{\la}^c$. Consider the isometric immersion of the form
$$\begin{array}{cccl}i_{\la}:=i_{\la}^0:&\Di{S}^2&\hookrightarrow&\Di{S}^4_1\subset\R^5_1.\\
&y&\mapsto&(\la,\la,y),\quad \|y\|=1\end{array}$$
This surface in $\Di{S}^4_1$ is the image of the surface $i_{\tilde{\la}}(\R^2)$ in $\R^4_1$ with $$\quad\tilde{\la}:=\exp(-\rho)\cdot\la=\frac{1+z_1^2+z_2^2}{2}\cdot\la$$ under the inverse of the stereographic projection $\Xi$:
$$\begin{array}{cccc}\Xi^{-1}\circ i_{\tilde{\la}}:&\R^2&\hookrightarrow&\Di{S}^4_1.\\[1.5mm]
&(z_1,z_2)&\mapsto&\frac{2}{1+z_1^2+z_2^2}(\tilde{\la},\tilde{\la},
\frac{z_1^2+z_2^2-1}{2},z_1,z_2)\\[2mm]
&&&=(\la,\la,\frac{z_1^2+z_2^2-1}{1+z_1^2+z_2^2},\frac{2z_1}{1+z_1^2+z_2^2},
\frac{2z_2}{1+z_1^2+z_2^2})\end{array}$$
We use this to calculate the second fundamental form $\tilde{II}$ of the isometric immersion
$i_{\la}:\Di{S}^2\hookrightarrow\Di{S}^4_1$. It holds (comp. Remark \ref{R3.1})
\begin{eqnarray*}\tilde{II}&=&<,>^{\Di{S}^2}\otimes[\exp(-2\rho)(H-
\Cal{N}grad(\rho))]\ +\ L_+\\
&=&<,>^{\Di{S}^2}\otimes[\De^{\Di{S}^2}\la+2\la]\cdot\frac{1}{2}\left(\frac{\partial}{\partial x_1}+\frac{\partial}{\partial x_2}\right)\ +\ L_+.\end{eqnarray*}
This proves that $i_{\la}:\Di{S}^2\hookrightarrow\Di{S}^4_1$ is a $+$isotropic immersion into $\Di{S}^4_1$. The immersions of the form $i_{\la}^c,\ c\neq 0$, are never $+$isotropic. We have the
\begin{THEO}
Any complete $+$isotropic surface in the pseudosphere $\Di{S}^4_1$ is isometric to $\Di{S}^2$ and up to an isometry of $\Di{S}^4_1$ a $+$isotropic surface has (locally) the form
$$\begin{array}{cccl}i_{\la}:&\Di{S}^2&\hookrightarrow&\Di{S}^4_1,\\
&y&\mapsto&(\la,\la,y),\quad \|y\|=1\end{array}$$
where $\la$ is a smooth function on $\R^2$. 
\end{THEO}
\textbf{Proof:} We need only to prove the global statement. The surface $i_{\la}(\Di{S}^2)\subset\Di{S}^4_1$ is compact and complete. Converserly, any $+$isotropic surface, which isn't a deformation of the complete sphere $i(\Di{S}^2)\subset\Di{S}^4_1$ can not be complete\hfill$\Box$\\
\begin{RE}
\begin{enumerate}
\item If $\hat{\la}$ is a spherical function to the eigenvalue $-2$ of the Laplace operator $\De^{\Di{S}^2}$ on $\Di{S}^2$, the surface $i_{\hat{\la}}(\Di{S}^2)\subset\Di{S}^4_1$ is complete, stationary and semi-umbilic.
\item
It is not true that all $+$isotropic surfaces in $\R^4_1$ and $\Di{S}^4_1$ are identified under a stereographic projection. 
\end{enumerate}
\end{RE}
C. Positive semi-umbilic and $+$isotropic surfaces in the pseudohyperbolic space $\Di{H}^4_1$\\[2mm]
Consider the pseudohyperbolic 4-space 
$$\Di{H}^4_1:=\{x\in\R^5_2:\ <x,x>^5_2=-1\}\subset\R^5_2.$$
The image of the imbedding
$$\begin{array}{cccl}\io:&\Di{H}^3\subset\R^4_1&\hookrightarrow&
\Di{H}^4_1\subset\R^5_2\\&(y_1,y_2,y_3,y_4)&\mapsto&(y_1,0,y_2,y_3,y_4),
\quad y_1>0,\ \|y\|=-1\end{array}$$
of the 3-dimensional hyperbolic space $\Di{H}^3$ into $\Di{H}^4_1$ is a complete spacelike hypersurface in $\Di{H}^4_1$, but there are lightlike geodesics in $\Di{H}^4_1$ that don't intersect $\io(\Di{H}^3)\subset\Di{H}^4_1$. Up to an isometry of $\Di{H}^3$, 
which is an element of $\SO_o(3,1)$, any totally umbilic surface in $\Di{H}^3$ has one of the following two forms:
$$\begin{array}{cccc}j:&\Di{S}^2&\hookrightarrow&\Di{H}^3,\\
&(z_1,z_2,z_3)&\mapsto&(c,\sqrt{c^2-1}\cdot z_1,\sqrt{c^2-1}\cdot z_2,\sqrt{c^2-1}\cdot z_3),
\quad c\geq 1\end{array}$$
$$\begin{array}{cccc}i:&\Di{H}^2&\hookrightarrow&\Di{H}^3.\\
&(z_1,z_2,z_3)&\mapsto&(z_1,0,z_2,z_3)\end{array}$$
The deformations of the surfaces $j(\Di{S}^2)\subset\Di{H}^4_1$ and $i(\Di{H}^2)\subset\Di{H}^4_1$ along the positive oriented normal null geodesics, which start on $j(\Di{S}^2)$ resp. $i(\Di{H}^2)$, look as follows:
$$\begin{array}{cccl}j_{\la}:&\Di{S}^2&\hookrightarrow&\Di{H}^4_1\\
&z=(z_1,z_2,z_3)&\mapsto&(\la,\la\sqrt{c^2-1}+c,(\sqrt{c^2-1}+\la c)z),\quad c\geq 1,\end{array}$$
and
$$\begin{array}{cccc}i_{\la}:&\Di{H}^2&\hookrightarrow&\Di{H}^4_1,\\
&(z_1,z_2,z_3)&\mapsto&(z_1,\la,\la,z_2,z_3)\end{array}$$
where $\la$ is a smooth function resp. on $\Di{S}^2$ and $\Di{H}^2$. Any positive semi-umbilic surface in $\Di{H}^4_1$ is locally up to an isometry of $\Di{H}^4_1$ an immersion of the form $j_{\la}$ or $i_{\la}$. For the second fundamental form of an immersion of the form $i_{\la}$, we have
$$\tilde{II}=<,>^{\Di{H}^2}\otimes\left[\De^{\Di{H}^2}\la-2\la\right]\cdot
\frac{1}{2}\left(\frac{\partial}{\partial x_2}+\frac{\partial}{\partial x_3}\right)\ +\ L_+.$$ Hence, the immersions of the form $i_{\la}$ are isometric and $+$isotropic, whereas the immersions $j_{\la}$ are only semi-umbilic.
If $\hat{\la}$ is a solution of the equation $\De^{\Di{H}^2}\hat{\la}-2\hat{\la}=0$, the immersion $i_{\hat{\la}}$ is even stationary.
\begin{THEO}
Any complete $+$isotropic surface in $\Di{H}^4_1$ is isometric to $\Di{H}^2$ and up to an isometry of $\Di{H}^4_1$ a +isotropic surface has
(locally) the form
$$\begin{array}{cccl}i_{\la}:&\Di{H}^2&\hookrightarrow&\Di{H}^4_1,\\
&(y_1,y_2,y_3)&\mapsto&(y_1,\la,\la,y_2,y_3),\quad y_1>0,\ -y_1^2+y_2^2+y_3^2=-1\end{array}$$
where $\la$ is a smooth function on $\Di{H}^2$. 
\end{THEO}

The discussion of the $+$isotropic surfaces in the space forms $\R^4_1$, $\Di{S}^4_1$ and $\Di{H}^4_1$ proves:
\begin{THEO}
Let $M^4_1$ be an oriented Lorentzian 4-manifold with constant sectional curvature.
\begin{enumerate}
\item Any smooth deformation of a $+$isotropic surface $N^2$ in $M^4_1$ along the positive oriented normal null geodesics of $N^2$ remains $+$isotropic.
\item
Any $+$isotropic surface in $M^4_1$ is locally a smooth deformation of a totally geodesic surface $N^2$ in $M^4_1$ along the positive oriented normal null geodesics of $N^2$.
\end{enumerate}
\end{THEO}

\newpage\noindent\Large\textbf{References}\\[2mm]\normalsize
\begin{description}
\item[{[Baum81]}] Helga Baum.
\it Spin-Strukturen und Dirac-Operatoren uber pseudo-Riemannschen Mannigfaltigkeiten. \rm Nummer 41 in Teubner-Texte zur Mathemtik. Teubner, 1981.
\item[{[Bes87]}] Arthur L. Besse.
\it Einstein manifolds, \rm Band 10 der Reihe \it Ergebnisse der Mathematik und ihrer Grenzgebiete, 3.Folge. \rm Springer, 1987.
\item[{[Jac90]}] Howard Jacobowitz. \it
An introduction to CR structures, \rm volume 32 of the series \it
Mathematical Surveys and Monographs. \rm
Amer. Math. Soc., 1990
\item[{[O'N83]}] Barret O'Neill. \it
Semi-Riemannian geometry. \rm Pure and Applied Mathematics. Akademic Press, 1983
\item[{[Tra85]}] Andrzej Trautmann. \it Optical structures in relativistic theories. \rm Soci\'{e}t\'{e} Math\'{e}matique de France, Ast\'{e}risque, hors s\'{e}rie, 1985, p. 401-420.
\item[{[Nur96]}] Pawel Nurowski. \it
Optical geometries and related structures. \rm J. Geometry and Physics 18(1996), p. 335-348.
\item[{[Fri81]}] Thomas Friedrich. \it Self-duality of Riemannian manifolds and connections. \rm In Teubner-Text Bd. 34, Teubner-Verlag Leipzig, 1981.
\item[{[FG83]}] Th. Friedrich and R. Grunewald. \it
On Einstein metrics on the twistor space of a four-dimensional Riemannian manifold. \rm Math. Nachr. 123(1985), p. 55-60. 
\item[{[AHS78]}] M.F. Atiyah, N.J. Hitchin and I.M. Singer. \it
Self-duality in Four-dimensional Riemannian geometry. \rm Proc. R.S.London A362
(1978), 425-461.
\item[{[Hit81]}] N. Hitchin. \it Kaehlerian twistor spaces. \rm Proc. London Math. Soc., III. Ser., 43(1981), p. 133-150.
\item[{[FK82]}] Th. Friedrich, H. Kurke. \it Compact four-dimensional self-dual Einstein manifolds with positive scalar curvature. \rm Math. Nachr. 106(1982)
271-299.
\item[{[ES85]}] Eells J., Salamon S., \it Twistorial constructions of harmonic maps of surfaces into four-manifolds, \rm Ann. Scuola Norm. Sup. Pisa 12(1985), 589-640.
\item[{[JR90]}] Jensen G.R., Rigoli M., \it Neutral surfaces in Neutral Four-spaces, \rm Le Matematiche Vol. XLV(1990) -Fasc. II, pp. 407-443.
\item[{[Fri84]}] Friedrich Th., \it On surfaces in Four-Spaces, \rm
Ann. Glob. Analysis and Geometry Vol. 2, No.3(1984), pp. 257-287.
\item[{[Bry82]}] Bryant, R.L., \it Conformal and minimal immersions of compact surfaces into the 4-sphere, \rm Joun. Diff. Geom. 17(1982), 455-473.
\item[{[Bob98]}] Bobienski M., \rm private communication, Warszawa 1998.
\item[{[Lei97]}] Leitner F., \it The twistor space of a Lorentzian manifold, \rm SfB 288 Preprint No. 314, http://www-sfb288.math.tu-berlin.de  
\item[{[Fri97]}]  Th. Friedrich. \it On superminimal surfaces. \rm Archivum Mathematicum vol. 33(1997), pp. 41-56.
\end{description}
\end{sloppypar}
\end{document}